\documentstyle[12pt]{article}
\input{amssym.def}

\makeatletter
\@addtoreset{equation}{section}

\makeatother
\begin{document}

\def\g{\frak{g}}
\def\ga{\hat{\g}}
\newtheorem{Lemma}{Lemma}[section]
\newtheorem{Theorem}{Theorem}[section]
\newtheorem{Example}{Example}[section]
\newtheorem{Definition}{Definition}[section]
\newtheorem{Remark}{Remark}[section]
\newtheorem{Proposition}{Proposition}[section]
\newtheorem{Corollary}{Corollary}[section]
\def\h{\frak {h}}
\def\ha{\hat{\h}}
\def\n{\frak {n}}
\def\na{\hat{\n}}
\def\nb{\bar{\n}}
\def\C{\Bbb {C}}
\def\N{\Bbb {N}}
\def\Z{\Bbb {Z}}
\newcommand{\be}{\begin{equation}}
\newcommand{\ee}{\end{equation}}

\setcounter{section}{-1}
\hfill q-alg/9504024
\begin{center}\LARGE Combinatorial  constructions of modules for
infinite-dimensional Lie algebras, II. Parafermionic
space.\\\vspace{5mm} {\Large Galin Georgiev}
\end{center}\vspace{8mm}
\begin{abstract}\begin{em} The standard modules  for an affine Lie
algebra $\ga$ have natural subquotients called parafermionic spaces --
the underlying spaces for the so-called parafermionic conformal field
theories associated with $\ga.$

   We study the case $\ga =
\widehat{sl}(n+1,\C)$ for  any positive integral level $k
\geq 2.$ Generalizing the $\cal Z$-algebra approach of Lepowsky,
Wilson and Primc, we construct  a  combinatorial basis for the parafermionic
spaces in terms of colored partitions. The parts of these partitions
represent ''Fourier coefficients'' of generalized vertex operators
(parafermionic currents) and can be interpreted as statistically
interacting quasi-particles of color $i,\;1\leq i \leq n,$ and charge
$s,\; 1\leq s \leq k-1.$ From a combinatorial point of view, these
bases are essentially identical with the bases for level $k-1$
principal subspaces given in [GeI].  In the particular case of the
vacuum module, the character (string function) associated with our
basis is the formula of Kuniba, Nakanishi and Suzuki [KNS] conjectured
in a Bethe Ansatz layout.

   New combinatorial characters are established  for the whole standard
vacuum $\ga$-modules.
\end{em}\end{abstract}

\section{Introduction}

\subsection{}We begin with a short outline of some results
{}from [GeI] which will be needed here.

   Let $\g := sl(n+1,\C)$ with a triangular decomposition
$\g = \n_{-} \oplus \h \oplus \n_{+}$ and simple roots
$\alpha_{i},$ $1\leq i \leq  n,$ where the indices reflect the roots'
location on the Dynkin diagram (the considerations below can be
carried out for any simple Lie algebra $\g$ over $\C$ of type A-D-E).
Denote by $x_{\pm\alpha_{i}},$ $1\leq i \leq n,$ the Chevalley
generators of $\g.$ Let $Q :=
\sum_{i=1}^{n} \Z \alpha_{i},$  $P := \sum_{i=1}^{n} \Z \Lambda_{i}$
be the root and weight lattice respectively, where $\Lambda _{i},$
$1\leq i \leq n,$ are the fundamental weights of $g.$ For some formal
variable $t,$ consider the untwisted affinizations $\ga := \g \otimes
\C[t, t^{-1}] \oplus \C c,$ $ \nb := \n \oplus \C[t, t^{-1}]$ with a
scaling operator $D:= -td/dt.$ Let $\hat{\Lambda}_{i},$ $0\leq i
\leq n,$  be the fundamental weights of $\ga,$ so that
$\hat{\Lambda}_{i} = \hat{\Lambda}_{0} + \Lambda_{i}$ when $ 1\leq i
\leq n.$ For another formal variable $z,$ denote by $X_{\pm
\alpha_{i}}(z) := \sum_{m \in \Z}(x_{\pm \alpha_{i}}\otimes t^{m})
z^{-m -1}$ the vertex operators (bosonic currents) of charge one,
corresponding to simple roots and their negatives. Define also
currents of higher charge $r \in
\Z_{+}$:  $X_{\pm r\alpha_{i}}(z) := X_{\pm \alpha_{i}}(z)^{r}$ (cf.
[GeI]). The operator-valued (''Fourier )
coefficients  of $X_{r\alpha_{i}}(z)$ are called quasi-particles of
color $i$ and charge $r.$

   For a given positive integral level $k$ (the eigenvalue of $c$),
let $L(k \hat{\Lambda}_{0})$ be the vacuum standard $\ga$-module,
i.e., the highest weight integrable module with highest weight $ k
\hat{\Lambda}_{0}$ (and highest weight vector $v(k
\hat{\Lambda}_{0})).$ For the sake of simplicity, we shall consider in this
Introduction only vacuum highest weights, although the results
are proved  for a larger class of highest
weights. Denoting by $U(\cdot)$ universal enveloping
algebra, recall from [FS] that $W(k
\hat{\Lambda}_{0}):= U(\nb)\cdot v(k  \hat{\Lambda}_{0})$ is called
principal subspace of the standard module. We constructed in [GeI] a
basis for the principal subspace: It is built by quasi-particles of
colors $i,$ $1\leq i \leq n,$ and charges $s,$ $ 1\leq s \leq k,$
acting on the highest weight vector $ v(k \hat{\Lambda}_{0}).$ A
straightforward counting of the basis resulted in the character
formula announced by Feigin and Stoyanovsky [FS]:
\be \mbox{\rm Tr}\,q^{D} \left| \begin{array}[t]{l} \\
W(k\hat{\Lambda}_{0}) \end{array} \right.  =
\sum_{\begin{array}{c}{\scriptstyle p_{1}^{(1)},\ldots , p_{1}^{(k)}
\geq 0}\\{\scriptstyle \ldots \ldots
\ldots\ldots}\\{\scriptstyle p_{n}^{(1)},\ldots
, p_{n}^{(k)} \geq 0 }\end{array}}
\frac{q^{\frac{1}{2}\sum_{l,m =
1, \ldots, n}^{s,t = 1,\ldots , k}\;A_{lm} B^{st}
p_{l}^{(s)}p_{m}^{(t)}}}{\prod_{i=1}^{n} \prod_{s=1}^{k}
(q)_{p_{i}^{(s)}} }, \label{0.1}\ee
where $(A_{lm})_{l,m = 1}^{n}$ is
the Cartan matrix of $\g,$ $B^{st} :=
\mbox{min}\{s,t\},$ $1\leq s,t \leq k,$ and for $p \in \Z_{+},$ $(q)_{p} :=
(1-q)(1-q^{2})\cdots (1-q^{p}),$ $(q)_{0}:= 1$ (for the character of a
principal subspace with more general highest weight, see formula
(5.27) [GeI]). This formula was interpreted in the Introduction of
[GeI] as a character for the Fock space of $nk$ different free bosonic
quasi-particles ($n$ different colors and $k$ different charges) with
an additional two-particle interaction $A_{lm} B^{st}$ between a
quasi-particle of color $l$ and charge $s$ and another quasi-particle
of color $m$ and charge $t.$ This interaction  can be interpreted as a
statistical interaction in the sense of Haldane [H].

   As indicated in [GeI], one can generate in a similar fashion a
basis for the whole standard module, employing in addition
quasi-particles corresponding to the negative simple roots. We shall
see in Proposition \ref{Pro 0.1} at the end this Introduction that it
is sufficient to add only quasi-particles of charge $k$ corresponding
to the negative simple roots (we might also refer to the latter as to
{\it quasi-antiparticles} of charge $-k).$

\subsection{} One of our objectives in this paper is the  generalization
(to higher rank affine Lie algebras) of the vertex operator
construction of parafermionic spaces, given by Lepowsky and Primc [LP]
for $\ga =
\widehat{sl}(2,\C).$ Their construction was the untwisted version of
the ground-breaking works on $\cal Z$-algebras of Lepowsky and Wilson [LW].

   Consider the so-called $\ga \supset \ha$
coset subspace $L(k \hat{\Lambda}_{0})^{\ha^{+}}$ of
$L(k \hat{\Lambda}_{0})$ consisiting of all the $\ha^{+}:= \h
\otimes t\C[t]$-invariants, i.e., the vectors annihilated by
$\ha^{+}.$ The parafermionic space $L(k
\hat{\Lambda}_{0})^{\ha^{+}}_{kQ}$ is defined as the spase of
$kQ$-coinvariants in $L(k \hat{\Lambda}_{0})^{\ha^{+}},$ i.e.,
$L(k\hat{\Lambda}_{0})^{\ha^{+}}_{kQ} :=
L(k\hat{\Lambda}_{0})^{\ha^{+}}/(\rho(kQ) - 1)\cdot
L(k\hat{\Lambda}_{0})^{\ha^{+}},$ where $\rho$ is a natural action of the
abelian group $kQ \cong Q.$  There is a natural
projection $\pi_{L(k
\hat{\Lambda}_{0})^{\ha^{+}}}$ such  that $L(k
\hat{\Lambda}_{0})^{\ha^{+}} = \rho(kQ)\cdot \pi_{L(k
\hat{\Lambda}_{0})^{\ha^{+}}}\cdot W(k\hat{\Lambda}_{0})$ (cf.
(\ref{1.6}) and (\ref{2.302}) below) -- this explains why structural
results for the principal subspace are easily carried  to the $\ga
\supset \ha$ subspace (and henceforth  to the
parafermionic space itself). For example, the bosonic current
$X_{r\alpha_{i}}(z)$ has simply to be replaced by $ \pi_{L(k
\hat{\Lambda}_{0})^{\ha^{+}}}\cdot X_{r\alpha_{i}}(z)$ -- the
latter equals up to a nonzero constant and a power of $z$ the
familiar parafermionic current $\Psi_{r\alpha_{i}}(z)$ (called
generalized vertex operator in the general setting of  [DL]) which
by analogy is said to have a color $i$ and charge $r.$ Roughly
speaking, one could think of the parafermionic current $\Psi$ as
obtained from the bosonic current $X$ by factoring a {\em free}
bosonic field (cf.  (\ref{2.14})).  As expected, the quasi-particles
of color $i$ and charge $r$ in the parafermionic setting will be the
operator-valued (''Fourier'') coefficients of the current
$\Psi_{r\alpha_{i}}(z).$ Note by the way that the parafermionic
counterpart of the simple relation $X_{r\alpha_{i}}(z) :=
X_{\alpha_{i}}(z)^{r}$ is more sophisticated and involves the familiar
binomial correction terms:
\be \Psi_{r\alpha_{i}}(z) = \left(\prod_{\begin{array}{c}
{\scriptstyle l,p=1}
\\{\scriptstyle  l>p}
\end{array}}^{r} (z_{l}-z_{p})^{\frac{\langle
\alpha_{i} ,\alpha_{i} \rangle }{k}}\right)
\Psi_{\alpha_{i}} (z_{r}) \cdots
\Psi_{\alpha_{i}}(z_{1})\left|\begin{array}{l}
 \\ z_{r} = \cdots = z_{1} = z \end{array}\right.,
\label{0.2} \ee
where the binomial terms  are to be expanded in nonnegative
inegral powers of the second variable (before the expression is
restricted on the hyperplane $z_{r} = \cdots = z_{1} = z).$ There are
other novelties in the parafermionic picture: The most important one
is probably the new constraint
\be \Psi_{k\alpha_{i}}(z) = \mbox{const}\,\rho(k\alpha_{i}),\;
\mbox{const}\in \C^{\times},\; 1\leq i \leq n, \label{0.3}\ee
which implies that the maximal allowed quasi-particle charge is
reduced from $k$ in the principal subspace picture to $k-1$ in the
parafermionic picture. Another subtlety is that the scaling operator
$D$ has to be shifted in the parafermionic setting by $- D^{\ha},$
where $D^{\ha}$ is the rescaled 0th
mode of the Virasoro algebra associated with the vertex operator
algebra $U(\ha)\cdot v(k\hat{\Lambda}_{0}),$ $\ha :=
\h \otimes
\C[t, t^{-1}] \oplus \C c$ of $\h$ (cf. (\ref{2.201})).

   We build a quasi-particle basis for
$L(k\hat{\Lambda}_{0})^{\ha^{+}}_{kQ}$ in section 4 (the corresponding
basis for $L(k\hat{\Lambda}_{0})^{\ha^{+}}$ is then obtained by
multiplying with $\rho(k\alpha),\; \alpha \in Q).$  As expected, it is
essentially the same as the [GeI] quasi-particle basis for the level
$k-1$ principal subspace $W((k-1)\hat{\Lambda}_{0}).$ The only
difference is that the two-particle interaction between a
quasi-particle of color $l$ and charge $s$ and another quasi-particle
of color $m$ and charge $t$ has a ''parafermionic '' shift $-
A_{lm}\frac{st}{k}$ and is therefore given by
\be A_{lm}\mbox{min}\{s,t\} - A_{lm}\frac{st}{k} =
A_{lm} A^{(-1)}_{st},\; 1\leq s,t \leq k-1,\label{0.4}\ee where
$(A_{lm})_{l,m = 1}^{n}$ is again the Cartan matrix of $\g =
sl(n+1,\C)$ and $(A^{(-1)}_{st})_{s,t=1}^{k-1}$ is the inverse of the
Cartan matrix of $sl(k,\C).$ In other words, the character formula
associated with our basis is the one conjectured by Kuniba, Nakanishi
and Suzuki [KNS]:
\be  \mbox{Tr}\,q^{D- D^{\ha}}\left | \begin{array}[t]{l} \\
L(k\hat{\Lambda}_{0})^{\ha^{+}}_{kQ}
\end{array} \right. = \sum_{\begin{array}{c}{\scriptstyle
p_{1}^{(1)},\ldots , p_{1}^{(k-1)} \geq 0}\\{\scriptstyle \ldots
\ldots
\ldots\ldots}\\{\scriptstyle p_{n}^{(1)},\ldots
, p_{n}^{(k-1)} \geq 0 }\end{array}}
\frac{q^{\frac{1}{2}\sum_{l,m =
1, \ldots, n}^{s,t = 1,\ldots , k-1}\;A_{lm} A_{st}^{(-1)}
p_{l}^{(s)}p_{m}^{(t)}}}{\prod_{i=1}^{n} \prod_{s=1}^{k-1}
(q)_{p_{i}^{(s)}} }. \label{0.5}\ee Note that a dilogarithm proof of
this formula has already been announced by Kirillov [Kir], so it might
be appropriate to emphasize that our goal here is not simply to find
another proof of this beautiful formula, but to reveal the underlying
conceptual structure.  This in-depth approach, although painful, will
pay off as we shall see for example in the subsequent Proposition 0.1.
Note also that our arguments work for more general dominant highest
weight (cf.  (\ref{4.1}) below)-- the associated character formula is
given in (\ref{5.7}) Section 5.  After an appropriate normalization,
one immediately obtaines from the above formulas combinatorial
expressions for the corresponding string functions.

   Note that if we restrict our attention to the particular case $\g =
sl(2,\C),$ the above formula (\ref{0.5}) is the celebrated Lepowsky-Primc
character [LP].  The underlying basis in
[LP] is nevertheless very different from ours: Their construction employs only
quasi-particles of charge one, governed by the so-called ''difference
$2$ at distance $k-1$'' condition, while we work with quasi-particles of
charge $r,$ $1\leq r
\leq k-1,$ governed by a ''difference $2r$ at distance $1$''
condition. A straightforward generalization of the Lepowsky-Primc
parafermionic basis for $\g = sl(3,\C)$ was given  in [P].

\subsection{}We continue  with the higher level
generalization of the new character formula for the level one standard
module $L(\hat{\Lambda}_{0})$ given in Proposition 0.1 [GeI] (for
simplicity, we restrict ourselves to vacuum modules only).

   As explained already, from the [GeI] basis for the principal
subspace $W(k\hat{\Lambda}_{0})$ we easily obtain a basis for the $\ga
\supset \ha$ subspace $L(k\hat{\Lambda}_{0})^{\ha^{+}}.$ The
latter immediately implies a basis for the whole standard module
$L(k\hat{\Lambda}_{0})
\cong U(\ha)\cdot v(k\hat{\Lambda}_{0}) \otimes
L(k\hat{\Lambda}_{0})^{\ha^{+}}$: The associated
character formula is (cf. (\ref{0.1})):
\be \mbox{\rm
Tr}\,q^{D}\left|
\begin{array}[t]{l} \\
L(k\hat{\Lambda}_{0}) \end{array} \right. = \frac{1}{(q)_{\infty}^{n}}
\sum_{\begin{array}{c}{\scriptstyle p_{1}^{(1)},\ldots , p_{1}^{(k-1)}
\geq 0}\\{\scriptstyle \ldots \ldots
\ldots\ldots}\\{\scriptstyle p_{n}^{(1)},\ldots
, p_{n}^{(k-1)} \geq 0 }\end{array}}
\frac{q^{\frac{1}{2}\sum_{l,m =
1, \ldots, n}^{s,t = 1,\ldots , k-1}\;A_{lm} B^{st}
p_{l}^{(s)}p_{m}^{(t)}}}{\prod_{i=1}^{n} \prod_{s=1}^{k-1}
(q)_{p_{i}^{(s)}} }
\cdot \sum_{\alpha \in Q} q^{\frac{k}{2}\langle \alpha , \alpha\rangle
+ \langle \alpha, \sum_{i=1}^{n} r_{i}\alpha_{i}\rangle},\label{0.6}\ee
where $r_{i} := \sum_{s=1}^{k-1}s p_{i}^{(s)}$ (cf.  (\ref{5.16}) for
a more general highest weight). Note that the sum over $Q$
incorporates the contribution of the factors $\rho(k\alpha),$ $\alpha
\in Q,$ and is easily expressed in terms of the classical theta
function of degree $k$:
\be  \sum_{\alpha \in Q} q^{\frac{k}{2}\langle \alpha , \alpha\rangle
+ \langle \alpha, \mu\rangle} = q^{-\frac{\langle \mu , \mu \rangle}{2k}}
\Theta_{\mu}(q)= q^{-\frac{\langle \mu , \mu \rangle}{2k}}
\sum_{\gamma \in Q + \frac{\mu}{k}} q^{\frac{k}{2}\langle
\gamma , \gamma \rangle}.\label{0.7}\ee
Although the presence of the theta function is somewhat intimidating,
the above expression  can  still be interpreted as
a trace along a combinatorial basis: This would be a semiinfinite
monomial basis built up from quasi-particles corresponding to
(positive) simple roots; cf. [FS] where the case  $\g = sl(2,\C)$ is
discussed.

   But there is yet another -- and probably the most natural -- way to
generate a basis for the whole standard module starting from the
principal subspace: Simply add quasi-particles $x_{-k\alpha_{i}}(m)$
corresponding to negative simple roots (we shall call
those anti-quasiparticles of charge $-k$) and take into account the
identity
\be (z_{2} - z_{1})^{r \langle \alpha_{i}, \alpha_{i} \rangle}
X_{-k\alpha_{i}}(z_{2})X_{r\alpha_{i}}(z_{1}) \left|
\begin{array}[t]{l} \\ z_{1} = z_{2} = z \end{array} \right. =
\mbox{const}\, X_{-(k-r)\alpha_{i}}(z),\;\;\mbox{const} \in
\C^{\times}, \label{0.8}\ee
for every simple root $\alpha_{i},\; 1\leq i \leq n,$ and charge $r,\;
1\leq r \leq k.$ In particular, this identity implies that all the
quasi-antiparticles of charge $-r,\; 1\leq r < k,$ can be generated by
anti-quasiparticles of charge $-k$ and usual quasi-particles.
Moreover, when $r=k,$ one gets an important new constraint between
quasiparticles of charge $k$ and anti-quasiparticles of charge $-k$
(the right-hand side is just a constant). Following the layout of
[GeI], it is now not difficult to generate a basis for the whole
standard module and write down the character formula associated with
it:

\begin{Proposition} One has the following $q$-character for the vacuum
standard module at level $k$:

\be \mbox{\rm Tr}\,q^{D} \left| \begin{array}[t]{l} \\
L(k\hat{\Lambda}_{0}) \end{array} \right.  =
\sum_{\begin{array}{c}{\scriptstyle p_{\pm 1}^{(1)},\ldots , p_{\pm 1}^{(k)}
\geq 0}\\{\scriptstyle \ldots \ldots
\ldots\ldots}\\{\scriptstyle p_{\pm n}^{(1)},\ldots
, p_{\pm n}^{(k)} \geq 0}\\{\scriptstyle p_{-i}^{(s)} = 0\;\forall s <
k,\; \forall i}\end{array}}
\frac{q^{\frac{1}{2}\sum_{l,m =
1, \ldots, n}^{s,t = 1,\ldots , k}\;A_{lm} B^{st} (p_{+l}^{(s)}-
p_{-l}^{(s)})(p_{+m}^{(t)}- p_{-m}^{(t)}) +
\sum_{l=1}^{n}p_{+l}^{(k)}p_{-l}^{(k)}  }}{\prod_{i=1}^{n}
\prod_{s=1}^{k} (q)_{p_{+i}^{(s)}}(q)_{p_{-i}^{(s)}} }, \label{0.9}\ee
where $(A_{lm})_{l,m = 1}^{n}$ is the Cartan matrix of $\g,$ $B^{st}
:=
\mbox{min}\{s,t\},$ $1\leq s,t \leq k,$ and for $p \in \Z_{+},$ $(q)_{p} :=
(1-q)(1-q^{2})\cdots (1-q^{p}),$ $(q)_{0}:= 1.$\label{Pro
0.1}\end{Proposition}\vspace{5mm}
In complete analogy with the level
one particular case (cf. [GeI] Proposition 0.1), the above formula
follows from (\ref{0.6}) and the ''Durfee rectangle'' combinatorial
identity [A]
\be \frac{1}{(q)_{\infty}} := \prod_{l \geq 0} (1- q^{l})^{-1} =
\sum_{\begin{array}[t]{c}\scriptstyle a, b \geq 0\\\scriptstyle a - b
= {\rm const}\end{array} }\frac{q^{ab}}{(q)_{a} (q)_{b}}.
\label{0.10}\ee
The basis underlying the above expression will be discussed in details
somewhere else.

   In the simplest
particular case $n = k =1,$ the right-hand side of the above character
reduces to
$$ \sum_{p_{+1}, p_{-1} \geq 0} \frac{q^{p_{+1}^{2} +
p_{-1}^{2} - p_{+1}p_{-1}}}{(q)_{p_{+1}}(q)_{p_{-1}}},$$
a formula which appeared first in [FS].

\subsection{} Since the above Proposition 0.1 is the higher level
generalization of Proposition 0.1 [GeI], the natural question arises,
what is the higher level generalization of Proposition 0.2 [GeI],
i.e., how to build a basis from intertwining operators corresponding
to fundamental weights and their negatives (this question
was posed long time ago by J.  Lepowsky)? This remains a difficult
open problem although a promising breakthrough in the simplest
particular case $\ga =
\widehat{sl}(2,\C)$ was recently made by Bouwknegt, Ludwig and
Schoutens [BLS1] (cf. also [BPS], [I]): Using explicitly the easy to
describe fusion algebra
of the fusion category for $\widehat{sl}(2,\C)$-modules, they built a
basis of the above type (the so-called spinon basis) and derived the
corresponding combinatorial character formula. As expected, the [BLS1]
character  is much more complicated than the particular
$\widehat{sl}(2,\C)$-case of formula (\ref{0.9}). The reason is
that (\ref{0.9}) reflects a basis built up from usual (as opposed to
intertwining!) vertex operators, namely, the ones
corresponding to simple roots. Of course, the latter do not interchange
different modules and for that matter, no knowledge of the fusion
algebra is needed.

   Let us note that the [BLS1] basis in principal gradation was
identified with the basis proposed in [FIJKMY] (cf.  [BLS2]).
Moreover, the [BLS1] character formula was later independently derived
using crystaline spinon basis, i.e., spinon basis for standard modules
over $U_{q}(\widehat{sl}(2,\C))$ at $q=0$ (cf. [NY], [ANOT]). No
attempts have been made yet to crystalize the ($q$-deformation of the)
basis underlying the above  character formula (\ref{0.9}).

   Since the character (\ref{0.5}) was originally conjectured by
Kuniba, Nakanishi and Suzuki [KNS] from Thermodynamic Bethe
Ansatz (TBA) considerations, it remains an open problem to understand
the connection between our vertex operator construction and TBA. One
possible approach, adopted by Melzer [M], Foda and Warnaar [FW], [W]
in the context of Virasoro algebra modules, is through
``finitization'' of the above $q$-series, i.e., representing them  as a
limit of certain $q$-polynomials.

   Finally, from the point of view of our approach, one of the most
important and challenging unsolved problems is of course to find a
generalization for arbitrary (nonintegral) levels and highest weights
and also, to build vertex operator bases for other coset spaces. This
will in particular provide vertex operator constructions for many
other modules over $W$-algebras. It will also cast some light upon the
connection with the approach of Berkovich and McCoy [BM] to Virasoro
algebra modules from the minimal series and to branching functions
associated with various other coset spaces considered in [KM], [KKMM],
[DKKMM], [KMM], [WP], [BG].

\subsection{}In Section 1 we define the
$\ga \supset \ha$ coset subspace of a standard $\ga$-module $(g =
sl(n+1,\C))$ at any positive integral level $k \geq 2.$ The
parafermionic space is defined as its natural quotientspace.  Section
2 introduces relative vertex operators (parafermionic currents).
Section 3 explains how the notions of quasi-particle and
quasi-particle momomial (from part I) have to be modified in the
current setting. In Section 4 we build a quasi-particle monomial basis
for the parafermionic space. Section 5 presents the corresponding
character formulas for the parafermionic space (string function) and
for the whole standard module. Tables 1 and 2 in the Appendix
illustrate Examples 4.1 and 4.2, Section 4. \vspace{10mm}

{\em Acknowledgments} We are indebted to James Lepowsky whose advices,
expertise and continuous encouragement made this work possible. Many
thanks are due to Barry McCoy for the intense  discussions,
the constant support and interest toward this work. We would also like
to thank Alexander Zamolodchikov.

\section{The $\ga  \supset  \ha $
 coset subspace and its  parafermionic quotientspace}

\hspace{3ex}Although the exposition below is selfcontained, we shall
follow the framework of [GeI] and use many of its notations, definitions and
results. We continue working with $\ga =
\widehat{sl}(n+1,\C)$ and a level $k$ dominant integral highest weight
$\hat{\Lambda} = k_{0}\hat{\Lambda}_{0} + k_{j}\hat{\Lambda}_{j}$ for
some $j,\; 1\leq j \leq n,$ and $k_{0}, k_{j} \in \N,\; k_{0} + k_{j}
= k \geq 2,$ where $\{\hat{\Lambda}_{l}\}_{l=0}^{n}$ are the
fundamental weights of $\ga.$ In other words, $\hat{\Lambda} = k
\hat{\Lambda}_{0} + \Lambda,$ where $\Lambda = k_{j} \Lambda_{j} \in P$
and $\{\Lambda_{l}\}_{l=1}^{n}$ are the fundamental weights of $\g.$
(All the objects introduced in Sections 1,2 and 3 have obvious
generalizations for any dominant integral highest weight.)

   Recall that $\ha := \h \otimes \C [t^{-1}, t] \oplus \C c \subset
\ga$ is the affinization of the Cartan
subalgebra $\h \subset \g.$ The {\em $\ga
\supset  \ha $ coset subspace} $L(\hat{\Lambda})^{\ha^{+}}$ of the
standard $\ga$-module $L(\hat{\Lambda})$ is
defined as the vacuum subspace for $\ha^{+},$ i.e., the linear span of
all the vectors $\{v \in L(\hat{\Lambda})| \ha^{+}\cdot v = 0\},$
where $\ha^{+} := \h \otimes t\C [t].$ A concrete realization of this
space can be given in terms of the vertex operator construction of
basic modules: Set $M(k):= U(\ha)\otimes_{U(\h \otimes \C [t] \oplus
\C c)} \C,$ with $\h \otimes
\C [t]$ acting trivially on $\C$ and $c$ acting as $k$ (cf. [GeI]
Preliminaries). Let $V_{P} := M(1) \otimes \C[P]$ and recall that by
the classical interpretation of $V_{P}$ as a $\ga$-module [FK], [S],
one has $V_{P} \cong \oplus _{j=0}^{n} L(\hat{\Lambda}_{j})$ (cf.
[GeI], Section 2). Therefore, the standard module $L(\hat{\Lambda})$
can be explicitly generated by the universal eneveloping algebra
$U(\ga)$ acting on the highest weight vector
$v(\hat{\Lambda})
\subset V_{P}^{\otimes k}$ through  the
$(k-1)$-iterate $\Delta^{k-1}$ of the standard coproduct $\Delta$ (cf.
[GeI], Sections 2, 5; no confusion can arise from the fact that we
denote by the same letter $\Delta$ the root system of $\g).$
Then the $\ga \supset \ha $ coset subspace  of $L(\hat{\Lambda})$ is of
course
\be  L(\hat{\Lambda})^{\ha^{+}} :=
\mbox{\rm span}_{\C} \{v \in L(\hat{\Lambda}) = U(\ga)\cdot
v(\hat{\Lambda}) | \ha^{+} \cdot v = 0\} \subset V_{P}^{\otimes
k}. \label{1.3}\ee
Note how easy it is to reconstruct the whole space $L(\hat{\Lambda})$ from
$L(\hat{\Lambda})^{\ha^{+}} $ because of the canonical isomorphisms of
$D$-graded linear spaces  [LW]
\begin{eqnarray}\label{1.4}U(\ha^{-}) \,\otimes \,
L(\hat{\Lambda})^{\ha^{+}} &\stackrel{\sim}{\rightarrow}
&L(\hat{\Lambda} ) \\ h \,\otimes \,u  &\mapsto &
h\cdot u, \nonumber \end{eqnarray}
and
\be S(\ha^{-})\cong  U(\ha^{-}) \cong M(k),\label{1.5}\ee
where $\ha^{-} := \ha \otimes
t^{-1}\C[t^{-1}], $ $S(\cdot )$ is a  symmetric algebra and $M(k)$ was
defined above. One can therefore  consider the projection
\be \pi_{L(\hat{\Lambda})^{\ha^{+}}}:  L(\hat{\Lambda}) \rightarrow
L(\hat{\Lambda})^{\ha^{+}},\label{1.6}\ee
given by the corresponding direct sum decomposition
\be  L(\hat{\Lambda}) =
L(\hat{\Lambda})^{\ha^{+}}  \oplus
 \ha^{-}U(\ha^{-}) \cdot L(\hat{\Lambda})^{\ha^{+}}. \label{1.7}
\ee
We can further reduce the $\ga \supset \ha$ coset subspace using its
natural structure of module for the root lattice $Q$:  Observe
that the map $\alpha \mapsto e^{\alpha},$ $\alpha \in Q,$ defines an
action of $Q$ on
$V_{P} = M(1) \otimes \C[P]$ (cf. [GeI], Section 2) by
restriction to  the
right factor. It thus commutes with the action of $\ha^{+}$ on $V_{P},$
which affects  only the left factor. Define a diagonal  action of the
sublattice $kQ \subset  Q$ $(kQ \cong Q)$ on
$V_{P}^{\otimes k}:$
\be k\alpha \mapsto \rho (k\alpha) := \underbrace{e^{\alpha}\otimes \ldots
\otimes e^{\alpha}}_{k\;factors}, \; \alpha \in Q. \label{1.8}\ee
Note that it   commutes with the action $\Delta^{k-1}(\ha^{+})$ of
$\ha^{+}$
and hence preserves  $L(\hat{\Lambda})^{\ha^{+}}
\subset V_{P}^{\otimes k}.$ Following [DL], Ch. 4, we consider the space of
$kQ$-coinvariants in  the $kQ$-module
$L(\hat{\Lambda})^{\ha^{+}}:$
\be L(\hat{\Lambda})^{\ha^{+}}_{kQ} :=
L(\hat{\Lambda})^{\ha^{+}}\left/\begin{array}{c} \\
\mbox{Span}\{(\rho(k\alpha) -1)\cdot v | \alpha \in Q,\; v \in
L(\hat{\Lambda})^{\ha^{+}}\} \end{array}\right. . \label{1.9}\ee This
quotientspace is called {\em parafermionic space} (of highest weight
$\hat{\Lambda}$) because it is a building block for the so-called
parafermionic conformal field theories [ZF], [G], [LP], [DL]. We shall
denote by
\be \pi_{L(\hat{\Lambda})^{\ha^{+}}_{kQ}}\, : \, L(\hat{\Lambda}) \rightarrow
L(\hat{\Lambda})^{\ha^{+}}_{kQ}\label{1.101}\ee
the composition of $\pi_{L(\hat{\Lambda})^{\ha^{+}}}$ from (\ref{1.6}) and
the obvious projection of $L(\hat{\Lambda})^{\ha^{+}}$ onto
$L(\hat{\Lambda})^{\ha^{+}}_{kQ}.$

   Our goal here is to construct a quasi-particle basis for the
parafermionic space $L(\hat{\Lambda})^{\ha^{+}}_{kQ}, $ modifying
appropriately the  quasi-particle basis for the principal subspace
$W(\hat{\Lambda}) \subset L(\hat{\Lambda})$ built in [GeI]. Observe
that if we denote by
$L_{\mu}(\hat{\Lambda})^{\ha^{+}}$ the weight subspace of
$L(\hat{\Lambda})^{\ha^{+}},$ corresponding to $\mu \in P$ and by $
L_{\mu}(\hat{\Lambda})^{\ha^{+}}_{kQ}$ its isomorphic (as
graded linear space) image in the quotient
$L(\hat{\Lambda})^{\ha^{+}}_{kQ},$ we have
\be \rho(k\alpha) \cdot L_{\mu}(\hat{\Lambda})^{\ha^{+}} =
L_{\mu+ k\alpha}(\hat{\Lambda})^{\ha^{+}},\;\alpha \in Q,\;\mu \in
P,\label{1.10}\ee
and
\be L(\hat{\Lambda})^{\ha^{+}}_{kQ} = \coprod_{\mu \in \Lambda +
Q/kQ} L_{\mu}(\Lambda)^{\ha^{+}}_{kQ} \cong  \coprod_{\mu
\in \Lambda +
Q/kQ} L_{\mu}(\hat{\Lambda})^{\ha^{+}}.\label{1.11}\ee
Therefore a basis for $L(\hat{\Lambda})^{\ha^{+}}_{kQ}$ furnishes
automatically a basis for the
whole $\ga \supset \ha$
coset subspace $L(\hat{\Lambda})^{\ha^{+}} = \coprod_{\mu \in \Lambda + Q}
L_{\mu}(\hat{\Lambda})^{\ha^{+}}.$ (No need to say,
the tensor product decomposition (\ref{1.4}) then provides  a basis
for the whole standard module $L(\hat{\Lambda}).)$

   The trivial action $\rho$ of $kQ$ on the parafermionic
space  is very suggestive for
considering a natural  action (denoted by the same letter $\rho)$ of
the finite abelian group $Q/kQ$ on
$L(\hat{\Lambda})^{\ha^{+}}_{kQ}$ (cf. [DL], Ch. 6): Set
\be \rho (\alpha) \cdot v := e^{2\pi i \frac{\langle \alpha , \mu
\rangle}{k}} v, \;\alpha \in Q/kQ,\; v \in
L_{\mu}(\hat{\Lambda})^{\ha^{+}}_{kQ}. \label{1.12}\ee
The characters $e^{2\pi i \frac{\langle \cdot , \mu
\rangle}{k}},\; \mu \in \Lambda + Q/kQ,$ are indeed the simple
characters of the group $Q/kQ,$ in other words, the decomposition
(\ref{1.11}) coincides with the character-space decomposition of
$L(\hat{\Lambda})^{\ha^{+}}_{kQ}$ under the above action of $Q/kQ.$

\section {Generalized  vertex operators (parafermionic currents)}
\hspace{3ex}In this Section, we shall largely use the methods of
generalized vertex operator algebra theory developed by Dong and
Lepowsky [DL] (cf. also [FLM]).

   Recall that the main protagonist in the level $k$ setting
of  [GeI] was the vertex
operator (bosonic current)
\be X_{\beta}(z) := \Delta^{k-1}(Y(e^{\beta},z)) = \label{2.1}\ee
$$= \underbrace{Y(e^{\beta}, z)\otimes \mbox{\bf 1}\otimes  \cdots
\otimes\mbox{\bf 1}}_{k\,factors} + \underbrace{\mbox{\bf 1}\otimes
Y(e^{\beta}, z) \otimes  \cdots
\otimes\mbox{\bf 1}}_{k\,factors} + \cdots + \underbrace{\mbox{\bf
1}\otimes \cdots
\otimes\mbox{\bf 1}\otimes Y(e^{\beta}, z)}_{k\,factors}, $$
where $\beta \in \Delta,$
\be Y(e^{\beta},z): =
E^{-}(-h_{\beta},z)E^{+}(-h_{\beta},z)\otimes
e^{\beta} z^{h_{\beta}} \varepsilon_{\beta},\label{2.101} \ee
\be
E^{\pm}(h,z) := \mbox{exp}\left( \sum_{m \geq 1} h(\pm m) \frac{z^{\mp
m}}{\pm m} \right), \; h \in \h, \label{2.2}\ee
and $\varepsilon_{\beta} := \varepsilon (\beta , \cdot ),$
$\varepsilon : P \times  P \rightarrow \C^{\times}$ being a 2-cocycle
on the weight lattice $P$ (cf. [GeI] Sections 1 - 3). The ''Fourier
coefficients'' of this vertex operator are given by its  expansion
\be X_{\beta}(z)=: \sum_{m \in \Z}x_{\beta}(m) z^{-m
-1},\label{2.102}\ee
on powers of the formal variable $z.$ The action of the affine
algebra $\ga$ on the standard module $L(\hat{\Lambda})
= U(\ga)\cdot v(\hat{\Lambda}) $ is then given by  $ x_{\beta}\otimes t^{m} :=
x_{\beta}(m),$ $\beta \in \Delta,\; m \in \Z.$

   The Jacobi
identity for the vertex operator algebra $L(k\hat{\Lambda}_{0})$
generated by these currents,  implies the  usual formulas for the
currents of higher charge  $r \in \Z_{+}:$
\be X_{r\beta}(z):= \underbrace{X_{\beta}(z) \cdots X_{\beta}(z)}_{r
factors} := Y( x_{\beta}(-1)^{r}\cdot
v(k\hat{\Lambda}_{0}),z),\label{2.3}\ee where $v(k\hat{\Lambda}_{0})$
is the vacuum highest weight vector at level $k$ (cf. for example [DL]
Proposition 13.16). Note that since the product $X_{\beta}(z_{2})
X_{\beta}(z_{1})$ is not singular on the hyperplane $z_{2} = z_{1} =
z$ (cf. [GeI] (2.8)), one does not need to ``regularize'' it with
powers of $z_{2} - z_{1}$ before one sets $z_{2} = z_{1} = z.$ This is
not the case anymore in the parafermionic picture -- see (\ref{2.13})
below.  We actually showed in [GeI] that the currents corresponding to
the simple roots $\alpha_{i},\; 1\leq i
\leq n,$ are enough to build a basis of the principal subspace
$W(\hat{\Lambda}):$ their ''Fourier coefficients'' played the role of
quasi-particles and the basis was generated by quasi-particle
monomials (from an appropriate completion of the ordered product $U :=
U(\nb_{\alpha_{n}})\cdots U(\nb_{\alpha_{1}})$ acting on the highest
weight vector $v(\hat{\Lambda})).$

   Switching now our attention from the principal subspace
$W(\hat{\Lambda})$ to the parafermionic space
$L(\hat{\Lambda})^{\ha^{+}}_{kQ},$ we immediately observe that the
above operators do not even preserve the vacuum subspace
$L(\hat{\Lambda})^{\ha^{+}},$ let alone its quotientspace
$L(\hat{\Lambda})^{\ha^{+}}_{kQ}!$ Fortunately, one can perform a
small cosmetic operation and fix this problem (cf. [LP], [ZF], [G],
[DL]): For every $\beta
\in \Delta,$ replace
$X_{\beta}(z)$ by the {\em relative vertex operator}
(called also {\em parafermionic current}) on $V_{P}^{\otimes k}$
\be   \Psi_{\beta}( z)
:= \left ( \underbrace{E^{-}(\frac{
 h_{\beta}}{k} , z) \otimes \cdots \otimes E^{-}(\frac{
h_{\beta}}{k} , z)}_{k\;factors} \right )
X_{\beta}(z) \label{2.4}\ee
$$\left  (
\underbrace{z^{-\frac{
h_{\beta}}{k}}\varepsilon_{\beta}^{-\frac{1}{k}}\otimes \cdots \otimes
z^{-\frac{h_{\beta }}{k}}\varepsilon_{\beta}^{-\frac{1}{k}}
}_{k\;factors} \right )\left ( \underbrace{E^{+}(\frac{ h_{\beta}}{k}
, z) \otimes\cdots \otimes E^{+}(\frac{ h_{\beta}}{k} ,
z)}_{k\;factors} \right ).$$ Note that the coefficients of this
operator lie in an appropriate completion of $U(\ga),$ i.e., they are
infinite sums which are truncated when acting on modules from the
category $\cal{O}.$ We do not have to specify here exactly which root
of unity is to be taken in $\varepsilon_{\beta}^{-\frac{1}{k}}$ as
long as it is the same on all tensor slots. The first and the last
factor on right-hand side of the above definition, together with the
$k^{th}$ root of formula (2.7) [GeI], ensure that
\be \left[ \ha^{+},  \Psi_{\beta}( z)\right] =  \left[
\ha^{-},  \Psi_{\beta}( z)\right] = 0,\label{2.5}\ee
i.e., the relative vertex operator indeed preserves the vacuum space
$L(\hat{\Lambda})^{\ha^{+}}$ (recall that $\ha$ acts on
$V_{P}^{\otimes k}$ through the iterated coproduct $\Delta^{k-1}).$
But it is the term in the middle (not present in the $\cal
Z$-operators of Lepowsky and Primc [LP])
\be \underbrace{z^{-\frac{
h_{\beta}}{k}}\varepsilon_{\beta}^{-\frac{1}{k}}\otimes \cdots \otimes
z^{-\frac{h_{\beta }}{k}}\varepsilon_{\beta}^{-\frac{1}{k}}
}_{k\;factors},\label{2.6}\ee
which guarantees that
\be \left[ \rho(k\alpha) ,  \Psi_{\beta}(z)\right] = 0,\; \alpha,
\beta \in Q,\label{2.7}\ee
and therefore $\Psi$ is well defined on the parafermionic space
$L(\hat{\Lambda})^{\ha^{+}}_{kQ}.$ It will be clear from the context
whether $\Psi$ acts on $L(\hat{\Lambda})^{\ha^{+}}$ or
$L(\hat{\Lambda})^{\ha^{+}}_{kQ}.$ \vspace{5mm}

\noindent {\it Remark 2.1} We should note that the nonzero
numerical coset correction $
\varepsilon_{\beta}^{-\frac{1}{k}}\otimes \cdots \otimes
\varepsilon_{\beta}^{-\frac{1}{k}}$ is not present in the definition
adopted in [DL]. As a result, the relative vertex operator in [DL]  does not
commute with $\rho (kQ)$ but only with $\rho (2kQ)$ and one is forced
to consider a larger parafermionic space associated with the finite
group $Q/2kQ.$\vspace{5mm}

   The components (''Fourier coefficients'') of $\Psi$  will be
indexed in the very
same fashion as the coefficients of $X,$ but
due to the term  (\ref{2.6}),
their indices  will  typically be rational numbers  rather than integers:
For every $\beta \in \Delta,$ set
\be \Psi_{\beta}(z) \left | \begin{array}[t]{l} \\
L_{\mu}(\hat{\Lambda})^{\ha^{+}} \end{array} \right.  =: \sum_{m\in \Z +
\frac{\langle \beta , \mu \rangle}{k}}
\psi_{\beta}(m) z^{-m -1}\left | \begin{array}[t]{l} \\
L_{\mu}(\hat{\Lambda})^{\ha^{+}} \end{array} \right. , \label{2.8}\ee
where $L_{\mu}(\hat{\Lambda})^{\ha^{+}}$ is
the $\mu$-weight subspace of $L(\hat{\Lambda})^{\ha^{+}}, $ $\mu
\in P$ (cf. for example [DL] (6.52); note that the operator
$\psi_{\beta}(m)$ is
defined only on those $\mu$-weight subspaces for which $m\in \Z +
\frac{\langle \beta , \mu \rangle}{k}).$   This definition is designed   so
that the coefficients $\psi_{\beta}(m)$ can also be
thought of as operators on the parafermionic space
$L(\hat{\Lambda})^{\ha^{+}}_{kQ}:$
\be \left[ \rho(k\alpha) ,  \psi_{\beta}(m)\right] = 0,\; \alpha,
\beta \in Q.\label{2.9}\ee
We call $-m -\frac{1}{k}$ (rather than $-m)$ a {\it conformal energy}
of $\psi_{\beta}(m).$ This is because the Virasoro algebra generators
are shifted on a coset space: In our parafermionic setting, we have to
replace the grading operator $D$ of the vertex operator algebra
$L(k\hat{\Lambda}_{0})$ (cf. [GeI], Preliminaries) by $D - D^{\ha},$
where
\be \left[ D - D^{\ha},  \psi_{\beta}(m)\right]  = - (m
 +\frac{1}{k})\psi_{\beta}(m) ,\label{2.105}\ee cf. e.g. [DL], (6.42)
and (14.87). (Formula (\ref{2.105}) is the parafermionic counterpart
of the commutation relation $[D, x_{\beta}(m)] = - m x_{\beta}(m).)$
In other words,
\be D^{\ha} \left | \begin{array}[t]{l} \\
L_{\mu}(\hat{\Lambda})^{\ha^{+}} \end{array} \right.  = L_{0}^{\ha} -
\frac{\langle \Lambda , \Lambda \rangle}{2k}\left | \begin{array}[t]{l} \\
L_{\mu}(\hat{\Lambda})^{\ha^{+}} \end{array} \right. =  \frac{\langle
\mu , \mu \rangle}{2k}  -
\frac{\langle \Lambda , \Lambda \rangle}{2k}\left | \begin{array}[t]{l} \\
L_{\mu}(\hat{\Lambda})^{\ha^{+}} \end{array} \right., \label{2.201}\ee
where $L_{0}^{\ha}$ is the $0$th mode of the Virasoro algebra
associated with the vertex operator algebra $U(\ha)\cdot
v(k\hat{\Lambda}_{0}) \cong M(k)$ (cf. [DL] (14.52)
and [K] Remark 12.8 for example).

   Pivotal for our arguments will be the observation that the
''Fourier coefficients'' of $\Psi_{\beta}(m)$ can be tied with the
''Fourier coefficients'' of $X_{\beta}(z)$ in a very simple way: For
any given weight $\mu \in P,$ one has from the very definitions
(\ref{2.102}), (\ref{2.4}) and (\ref{2.8})
\be  \pi_{L(\hat{\Lambda})^{\ha^{+}}} \cdot x_{\beta}(m)\left |
\begin{array}[t]{l} \\
L_{\mu}(\hat{\Lambda})^{\ha^{+}} \end{array} \right.  =
\mbox{const}\,
\psi_{\beta} (m + \frac{\langle \beta , \mu \rangle}{k})\left |
\begin{array}[t]{l} \\
L_{\mu}(\hat{\Lambda})^{\ha^{+}} \end{array} \right., \label{2.10}\ee
where $m \in \Z,$ $\mbox{const} \in \C^{\times}$ and the natural projection
$\pi_{L(\hat{\Lambda})^{\ha^{+}}}: L(\hat{\Lambda}) \rightarrow
L(\hat{\Lambda})^{\ha^{+}}$ was introduced in (\ref{1.6}),
(\ref{1.7}). Note that this identity does not make much sense if
$L(\hat{\Lambda})^{\ha^{+}}$ is replaced by
$L(\hat{\Lambda})^{\ha^{+}}_{kQ}$ and if $\psi$ is thought of as an
operator on $L_{\mu}(\hat{\Lambda})^{\ha^{+}}_{kQ}$ because
$\pi_{L(\hat{\Lambda})^{\ha^{+}}_{kQ}} \cdot x_{\beta}(m)$ does not
really act on the quotientspace $L(\hat{\Lambda})^{\ha^{+}}_{kQ}.$
This is one of the reasons why we shall often use $L(\hat{\Lambda})^{\ha^{+}}$
as a mediator between  the principal subspace and the parafermionic space.

   Observe  that one can reverse the definition (\ref{2.4}) and thus
''isolate'' the part of the vertex
operator $X_{\beta}(z)$ which acts on the $\ga \supset \ha$ subspace
$L(\hat{\Lambda})^{\ha^{+}} \subset L(\hat{\Lambda}):$
\be X_{\beta}(z) =  \left ( \underbrace{E^{-}(-\frac{
 h_{\beta}}{k} , z) \otimes \cdots \otimes E^{-}(-\frac{
h_{\beta}}{k} , z)}_{k\;factors} \right ) \Psi_{\beta}( z) \label{2.11}\ee
$$\left  (
\underbrace{z^{\frac{
h_{\beta}}{k}}\varepsilon_{\beta}^{\frac{1}{k}}\otimes \cdots \otimes
z^{\frac{h_{\beta }}{k}}\varepsilon_{\beta}^{\frac{1}{k}}
}_{k\;factors} \right )\left (
\underbrace{E^{+}(-\frac{
 h_{\beta}}{k} , z) \otimes\cdots  \otimes E^{+}(-\frac{
h_{\beta}}{k} , z)}_{k\;factors} \right ).$$\vspace{5mm}
\begin{Lemma} The $\ga \supset \ha$ subspace
$L(\hat{\Lambda})^{\ha^{+}}$
(and therefore the parafermionic space $L(\hat{\Lambda})^{\ha^{+}}_{kQ})$ is
generated by  the operators $\{\psi_{\beta}(m) |
\beta \in \Delta \}$ acting on the highest weight vector $v(\hat{\Lambda}),$
i.e.,
\be L(\hat{\Lambda})^{\ha^{+}} =
\mbox{\rm Span}_{\C}\left\{\psi_{\beta_{r}}(m_{r}) \cdots
\psi_{\beta_{1}}(m_{1}) \cdot v(\hat{\Lambda}) | \beta_{l} \in \Delta, \;1
\leq l \leq r\right\}.\label{2.12}\ee\label{Lem 2.1}\end{Lemma}

\noindent{\em Proof}  Assume the opposite and using (\ref{2.11}) and
(\ref{1.4}),  arrive at a contradiction with the irreducibility of
$L(\hat{\Lambda})$ (cf. [DL] Proposition 14.9).\hfill $\Box$\vspace{5mm}

   Note by the way  that for $\beta \in \Delta,$
\be\rho(k\beta)\cdot v(\hat{\Lambda}) =  \mbox{const}\, x_{\beta}(-1
-\langle \beta , \Lambda\rangle) x_{\beta}(-1)
^{k-1} \cdot v(\hat{\Lambda}) = \label{2.202}\ee
$$= \mbox{const}'\,
\pi_{L(\hat{\Lambda})^{\ha^{+}}} \cdot x_{\beta}(-1
-\langle \beta , \Lambda\rangle) x_{\beta}(-1)
^{k-1} \cdot v(\hat{\Lambda}) = $$
$$=  \mbox{const}''\,
\pi_{L(\hat{\Lambda})^{\ha^{+}}} \cdot  x_{\beta}(-1
-\langle \beta , \Lambda\rangle) \left(
\pi_{L(\hat{\Lambda})^{\ha^{+}}}\cdot  x_{\beta}(-1)\right)^{k-1}
\cdot v(\hat{\Lambda})$$
for some nonzero constants. Since $v(\hat{\Lambda})$ is an eigenvector
for $D - D ^{\ha}$ (cf. (\ref{2.201})), one can compute in a
straightforward fashion from (\ref{2.105}), (\ref{2.10}) and
(\ref{2.202}) that $D - D^{\ha}$ and $\rho(k\alpha),\;
\alpha \in Q,$ commute when acting on a highest weight vector. It
follows immediately  from (\ref{2.9}), (\ref{2.105}) and Lemma
\ref{Lem 2.1} that $\rho(k\alpha)$ and $D - D^{\ha}$ commute on
$L(\hat{\Lambda})^{\ha^{+}}.$ But $\rho(k\alpha) $ acts nontrivially
only on the right factor of the decomposition (\ref{1.4}), hence
\be [D - D^{\ha},\rho(k\alpha)] = 0, \; \alpha \in Q.\label{2.203}\ee

   Let us continue now with the parafermionic currents of higher
charge (well known in the physics literature -- cf. e.g. [ZF], [G]): For a
given $\beta \in \Delta$ and $r \in \Z_{+},$ we call a {\em parafermionic
current of charge $r$} the generating function
\be \Psi_{r\beta}(z) := \left(\prod_{\begin{array}{c}
{\scriptstyle l,p=1}
\\{\scriptstyle  l>p}
\end{array}}^{r} (z_{l}-z_{p})^{\frac{\langle
\beta ,\beta \rangle }{k}}\right)
\Psi_{\beta} (z_{r}) \cdots \Psi_{\beta}(z_{1})\left|\begin{array}{l}
\\ z_{r} = \cdots = z_{1} = z \end{array}\right.,
\label{2.13} \ee
where the binomial terms  are to be expanded in nonnegative
inegral powers of the second variable (before the expression is
restricted on the hyperplane $z_{r} = \cdots = z_{1} = z).$ This
generating function is
well defined when acting on a highest weight module because the
binomial terms  cancel exactly the singularities related to  the
noncommutativity of the first and the last correction factors in the definition
(\ref{2.4}) of  $\Psi$ (cf. [GeI] (2.7),
(2.8)). In other words, the above expression can be rewritten as
\be  \Psi_{r\beta}(z) = \mbox{const} \left(
\underbrace{E^{-}(\frac{
 h_{\beta}}{k} , z) \otimes \cdots  \otimes \,E^{-}(\frac{
h_{\beta}}{k} , z)}_{k\;factors} \right )^{r} X_{r\beta}(z)\label{2.14}\ee
$$ \left(
\underbrace{z^{-\frac{
h_{\beta}}{k}}\varepsilon_{\beta}^{-\frac{1}{k}}\otimes \cdots\otimes
z^{-\frac{h_{\beta}}{k}}\varepsilon_{\beta}^{-\frac{1}{k}}}_{k\;factors}
\right)^{r} \left(
\underbrace{E^{+}(\frac{
 h_{\beta}}{k} , z) \otimes \cdots \otimes E^{+}(\frac{ h_{\beta}}{k}
, z)}_{k\;factors} \right )^{r},$$ where $\mbox{const}\in \C^{\times}$
and $ X_{r\beta}(z) = X_{\beta}(z)^{r}$ is the bosonic current of
charge $r$ from (\ref{2.3}). It will be clear from the context whether
$\Psi$ acts on $L(\hat{\Lambda})^{\ha^{+}}$ or
$L(\hat{\Lambda})^{\ha^{+}}_{kQ}.$

   The generating function  $\Psi_{r\beta}(z)$ is the   parafermionic
counterpart of $X_{r\beta}(z)$ in  the following sense:
Recall that $X_{r\beta}(z)$ is the vertex operator corresponding to
the vector  $x_{\beta}(-1)^{r}\cdot v(k
\hat{\Lambda}_{0} ) $ in the vertex operator algebra $ L(k
\hat{\Lambda}_{0} )$ (cf. (\ref{2.3})). According to (\ref{2.10}), the
projection of this vector on the parafermionic space is
\be \pi_{L(\hat{\Lambda})^{\ha^{+}}} \cdot  x_{\beta}(-1)
^{r}\cdot v(k
\hat{\Lambda}_{0} ) =  \mbox{const}\, \underbrace{\psi_{\beta}(-1 +
\frac{\langle \beta , (r-1)\beta\rangle}{k} )\cdots
\psi_{\beta}(-1)}_{r\, factors} \cdot v(k
\hat{\Lambda}_{0} ),\label{2.16}\ee
(the nonzero $\mbox{const}$ equals one if we   assume without losing
generality that
$\varepsilon$ is  bimultiplicative and $\varepsilon(\alpha , \alpha)
= \varepsilon(\alpha , \alpha)^{\frac{1}{k}} = 1).$  But the Jacobi
identity for the generalized (in the sense of [DL])
vertex operator algebra $L(k
\hat{\Lambda}_{0} )^{\ha^{+}}_{kQ}$ easily implies that
$\Psi_{r\beta}(z)$ is the vertex operator corresponding to this last
vector, i.e.,
\be \Psi_{r\beta}(z) = Y(\underbrace{\psi_{\beta}(-1 +
\frac{\langle \beta , (r-1)\beta\rangle}{k} )\cdots
\psi_{\beta}(-1)}_{r\,factors} \cdot v(k
\hat{\Lambda}_{0} ), z) \label{2.17}\ee
(one can for example use repeatedly  Proposition 14.29 [DL] which under
the above assumptions for $\varepsilon$ is still true
even for our slightly modified $\Psi,$ provided the notations are
appropriately adjusted;  be aware that the above
$Y(\cdot ,z)$ is not the same as $Y (\cdot ,z)$ in (\ref{2.3})
since these are vertex operators in two {\em different} vertex operator
algebras).

   In view of the correspondence between parafermionic and
bosonic currents, the most natural generalization of the definition
(\ref{2.8}) of ''Fourier coefficients'' for  higher-charge
parafermionic currents is (cf. [GeI] (3.7))
\be \Psi_{r\beta}(z) \left | \begin{array}[t]{l} \\
L_{\mu}(\hat{\Lambda})^{\ha^{+}} \end{array} \right.  =: \sum_{m\in \Z +
\frac{\langle r \beta , \mu \rangle}{k}}
\psi_{r\beta}(m) z^{-m -r} \left | \begin{array}[t]{l} \\
L_{\mu}(\hat{\Lambda})^{\ha^{+}} \end{array} \right. , \label{2.103}\ee
where $ \beta \in \Delta,\; r \in \Z_{+}$ and
$L_{\mu}(\hat{\Lambda})^{\ha^{+}}$ is  as always
the $\mu$-weight subspace of $L(\hat{\Lambda})^{\ha^{+}}, $ $\mu
\in P.$ Note that $\psi_{r\beta}(m)$ is a  well defined operator on
$L_{\mu}(\hat{\Lambda}
)^{\ha^{+}}_{kQ}$ with $\mu, $ such that $m\in \Z +
\frac{\langle r \beta , \mu \rangle}{k},$ because  $\psi_{r\beta}(m)$
commutes with the action $\rho$ of $kQ.$  Moreover, by (\ref{2.105})
and the definitions
(\ref{2.13}), (\ref{2.103}), the conformal energy of $\psi_{r\beta}(m)$ is
\be \left[ D - D^{\ha},  \psi_{r\beta}(m)\right]  = - (m
 + \frac{1}{k}(r + \langle \beta , \beta \rangle \left( \begin{array}{c}
 r\\2 \end{array}\right))\psi_{r\beta}(m) = \label{2.106}\ee $$= - (m  +
\frac{r^{2}}{k})\psi_{r\beta}(m).$$
Analogously to the charge-one situation (\ref{2.10}), one can find for
every operator $x_{r\beta}(m)$ its parafermionic counterpart (acting
on a given weight subspace of $L(\hat{\Lambda})^{\ha^{+}})$ by
composing it with the projection $\pi_{L(\hat{\Lambda})^{\ha^{+}}}:$
\be  \pi_{L(\hat{\Lambda})^{\ha^{+}}} \cdot x_{r\beta}(m)\left |
\begin{array}[t]{l} \\
L_{\mu}(\hat{\Lambda})^{\ha^{+}} \end{array} \right.  =
\mbox{const}\,
\psi_{r\beta} (m + \frac{ \langle r\beta , \mu \rangle}{k})\left |
\begin{array}[t]{l} \\
L_{\mu}(\hat{\Lambda})^{\ha^{+}} \end{array} \right.  ,\label{2.104}\ee
where $\mbox{const} \in \C^{\times}.$

   We are now in position to formulate a key relation which -- simple
and beautiful as it is -- explains both the similarity and the
difference between the structure of the parafermionic space
$L(\hat{\Lambda} )^{\ha^{+}}_{kQ}$ and the principal space
$W(\hat{\Lambda}).$\vspace{5mm}
\begin{Proposition} For every  $\beta \in \Delta$ and $r \in \Z_{+},\;
1\leq r \leq k,$
one has
\be \Psi_{r\beta}(z) =
\mbox{\rm const}\,\rho(k\beta) \Psi_{-(k -r)\beta} (z),\label{2.18}
\end{equation}
$ \mbox{\rm const} \in
\C^{\times},$ (cf. (\ref{1.8})), i.e.,
\be \Psi_{r\beta}(z) =
\mbox{\rm const}\, \Psi_{-(k -r)\beta} (z),\label{2.19} \end{equation}
$\mbox{\rm const} \in \C^{\times},$ as operators  on $L(\hat{\Lambda}
)^{\ha^{+}}_{kQ}.$ \label{Pro 2.2}\end{Proposition}

\noindent{\em Proof } Follows from  a direct computation employing
the commutation relation [GeI] (2.7) (in complete analogy with
[LP] Theorem 5.6 for example).\hfill $\Box$ \vspace{5mm}

   This in particular implies that products of coefficients of
parafermionic currents associated with positive roots are enough for
generating the parafermionic space (when acting on the vacuum vector;
cf. Lemma
\ref{Lem 2.1}
and (\ref{2.13})).
In other words, in complete  analogy with the principal (sub)space
\be W(\hat{\Lambda}) := \mbox{Span}_{\C}\left\{x_{\beta_{r}}(m_{r}) \cdots
x_{\beta_{1}}(m_{1}) \cdot v(\hat{\Lambda}) | \beta_{l} \in \Delta_{+}, \; 1
\leq l \leq r\right\},\label{2.20}\ee
one has
\be L(\hat{\Lambda} )^{\ha^{+}} = \rho(kQ)\cdot
\mbox{Span}_{\C}\left\{\psi_{\beta_{r}}(m_{r}) \cdots
\psi_{\beta_{1}}(m_{1}) \cdot v(\hat{\Lambda}) |  \beta_{l} \in
\Delta_{+}, \; 1
\leq l \leq r\right\}\label{2.301}\ee
and hence
\be L(\hat{\Lambda} )^{\ha^{+}}_{kQ} =
\mbox{Span}_{\C}\left\{\psi_{\beta_{r}}(m_{r}) \cdots
\psi_{\beta_{1}}(m_{1}) \cdot v(\hat{\Lambda}) |  \beta_{l} \in
\Delta_{+}, \; 1
\leq l \leq r\right\}.\label{2.21}\ee
Despite this similarity, there is an important difference  between the
parafermionic space and the
principal (sub)space which is encoded  in the particular case of
Proposition \ref{Pro 2.2} for $r = k$:
\be  \Psi_{k\beta}(z) = \mbox{const}\,
\rho(k\beta),\label{2.22}\ee
$\mbox{const} \in \C^{\times},\; \beta \in \Delta,$ i.e.,
\be  \Psi_{k\beta}(z) = \mbox{const} \in \C^{\times} \label{2.23}\ee
when acting on the parafermionic space $L(\hat{\Lambda}
)^{\ha^{+}}_{kQ}.$ In other words, one component of $\Psi_{k\beta}(z)$
acts as a
nonzero constant on the parafermionic space
$L(\hat{\Lambda} )^{\ha^{+}}_{kQ}$ and all the other components
vanish on it. Recall that on the principal subspace
$W(\hat{\Lambda}),$ we have the constraint $X_{(k+1)\beta}(z) = 0,$
which of course implies $\Psi_{(k+1)\beta}(z) = 0.$ In contrast to
this mutually shared constraint, the above constraint (\ref{2.23}) is
a purely parafermionic phenomenon with no analog in the principal
subspace. It tells us that the maximal allowed charge of the (defined
below) quasi-particles generating the parafermionic space is $k-1.$
This is also the maximal allowed charge of the quasi-particles
generating a principal subspace at level $k-1.$ So, not surprisingly,
the constructed below basis for a parafermionic space at level $k$
will be combinatorially the same as a [GeI] basis for a principal space
at level $k-1.$

   For the purpose of building a basis for $L(\hat{\Lambda}
)^{\ha^{+}}$ from the known already basis of $W(\hat{\Lambda}),$
it is very natural to employ not only the parafermionic counterpart
$\Psi_{r\beta}(z)$ of the bosonic current $X_{r\beta}(z),$ but also
the parafermionic counterpart of a whole product of bosonic currents
(with {\em different} variables) which differs from the product of the
respective parafermionic counterparts: For given roots $\beta_{r},
\ldots ,\beta_{1} \in \Delta$ and corresponding sequence of charges
$n_{r}, \ldots ,n_{1} \subset  \Z_{+},$ set
\be \Psi_{n_{r}\beta_{r}, \ldots , n_{1}\beta_{1}}(z_{r}, \ldots,
z_{1}) := \label{2.24}\ee
$$ := \left(\prod_{\begin{array}{c}
{\scriptstyle l,p=1}
\\{\scriptstyle  l>p}
\end{array}}^{r} (z_{l}-z_{p})^{\frac{\langle
n_{l}\beta_{l} ,n_{p}\beta_{p} \rangle }{k}}\right)
\Psi_{n_{r}\beta_{r}} (z_{r}) \cdots \Psi_{n_{1}\beta_{1}}(z_{1}),$$
where the binomial terms are to be expanded as usual in
nonnegative integral powers of the second variable. Just like in
(\ref{2.13}), they are inserted in order to ensure that the
composition $\pi_{L(\hat{\Lambda} )^{\ha^{+}}}\cdot
X_{n_{r}\beta_{r}}(z_{r})\cdots
X_{n_{1}\beta_{1}}(z_{1})$ equals (up a nonzero constant) $
\Psi_{n_{r}\beta_{r}, \ldots , n_{1}\beta_{1}}(z_{r}, \ldots
z_{1}).$ In other words, analogously to (\ref{2.14}), one can show
that
\be  \Psi_{n_{r}\beta_{r}, \ldots , n_{1}\beta_{1}}(z_{r}, \ldots
z_{1})  = \mbox{const} \prod_{l =1}^{r}\left(
\underbrace{E^{-}(\frac{
 h_{\beta_{l}}}{k} , z) \otimes \cdots  \otimes \,E^{-}(\frac{
h_{\beta_{l}}}{k} , z)}_{k\;factors} \right )^{n_{l}}\label{2.25}\ee
$$ X_{n_{r}\beta_{r}}(z_{r}) \cdots
X_{n_{1}\beta_{1}}(z_{1}) \prod_{l =1}^{r} \left(
\underbrace{z_{l}^{-\frac{
h_{\beta_{l}}}{k}}\varepsilon_{\beta_{l}}^{-\frac{1}{k}}\otimes \cdots\otimes
z_{l}^{-\frac{h_{\beta_{l}}}{k}}
\varepsilon_{\beta_{l}}^{-\frac{1}{k}}}_{k\;factors}
\right)^{n_{l}}$$
$$ \prod_{l =1}^{r}  \left(
\underbrace{E^{+}(\frac{
 h_{\beta_{l}}}{k} , z) \otimes  \cdots   \otimes  E^{+}(\frac{
h_{\beta_{l}}}{k} , z)}_{k\;factors} \right )^{n_{l}},$$
where $\mbox{const}\in \C^{\times}.$ We call $\Psi_{n_{r}\beta_{r},
\ldots , n_{1}\beta_{1}}(z_{r}, \ldots z_{1})$ a {\em normalized
product} of the parafermionic currents $ \Psi_{n_{r}\beta_{r}}
(z_{r}), \ldots , \Psi_{n_{1}\beta_{1}}(z_{1}).$ The last equality
guarantees that a normalized product is invariant (up to a nonzero
constant) under the permutation of two adjacent variables $z_{l}, z_{l+1}$
and the corresponding indices $n_{l}\beta_{l}, n_{l+1}\beta_{l+1}$ as
long as $\beta_{l} = \beta_{l+1}$ (this  is not true for the usual
product of parafermionic currents).  Note that for any $r \in \Z_{+},$
one has
\be \Psi_{\underbrace{\scriptstyle \beta, \ldots ,
\beta}_{r\,entries} }(z, \ldots , z)=
\Psi_{r\beta}(z).
\label{2.26}\ee
Generalizing (\ref{2.103}), one defines the ''Fourier coefficients''
of a normalized product as follows
\be \Psi_{n_{r}\beta_{r},
\ldots , n_{1}\beta_{1}}(z_{r}, \ldots z_{1}) \left | \begin{array}[t]{l} \\
L_{\mu}(\hat{\Lambda})^{\ha^{+}} \end{array} \right.  =: \label{2.27}\ee
$$=: \sum_{m_{r}\in \Z +
\frac{\langle n_{r} \beta_{r} , \mu \rangle}{k}}\cdots \sum_{m_{1}\in \Z +
\frac{\langle n_{1} \beta_{1} , \mu \rangle}{k}}
 \psi_{n_{r}\beta_{r},
\ldots , n_{1}\beta_{1}}(m_{r}, \ldots , m_{1})\left | \begin{array}[t]{l} \\
L_{\mu}(\hat{\Lambda})^{\ha^{+}} \end{array} \right.$$
$$ z_{r}^{-m_{r} -n_{r}}\cdots z_{1}^{-m_{1}-n_{1}} , $$
where $ \beta \in \Delta,\; r \in \Z_{+}$ and
$L_{\mu}(\hat{\Lambda})^{\ha^{+}}$ is  as always
the $\mu$-weight subspace of $L(\hat{\Lambda})^{\ha^{+}}, $ $\mu
\in P.$ We call $\psi_{n_{r}\beta_{r},
\ldots , n_{1}\beta_{1}}(m_{r}, \ldots , m_{1})$ a {\em normalized
$\psi$-monomial} or simply a {\em normalized
monomial}, as opposed to the usual ($\psi$-) monomial
$\psi_{n_{r}\beta_{r}}(m_{r}) \cdots \psi_{n_{1}\beta_{1}}(m_{1}).$
According to the very definition (\ref{2.24}), a normalized monomial
is a linear combination of usual monomials and vice versa (cf.
(\ref{3.11}), (\ref{3.12}) below; no need to
mention that our monomials are always acting on a designated space and
the linear combinations in question are truncated, i.e., finite). The
conformal energy of
a normalized  monomial is given by the following corollary of
(\ref{2.106}) and (\ref{2.24}):
\be \left[ D - D^{\ha}, \psi_{n_{r}\beta_{r},
\ldots , n_{1}\beta_{1}}(m_{r}, \ldots , m_{1}) \right]  =\label{2.28}\ee
$$ = -\sum_{l=1}^{r}(m_{l}  +  \frac{n_{l}^{2}}{k} +
\frac{\langle
n_{l}\beta_{l} , \sum_{p<l} n_{p}\beta_{p}}{k}) \psi_{n_{r}\beta_{r},
\ldots , n_{1}\beta_{1}}(m_{r}, \ldots , m_{1}) .$$
As expected from the preceding discussion, the normalized monomials
acting on $L(\hat{\Lambda})^{\ha^{+}}$
are indeed the parafermionic counterparts of monomials of type
$x_{n_{r}\beta_{r}}(m_{1})
\cdots x_{n_{1}\beta_{1}}(m_{1})$: One has from (\ref{2.25}) and
(\ref{2.27}) that
\be  \pi_{L(\hat{\Lambda})^{\ha^{+}}} \cdot x_{n_{r}\beta_{r}}(m_{r})
\cdots x_{n_{1}\beta_{1}}(m_{1}) \left |
\begin{array}[t]{l} \\
L_{\mu}(\hat{\Lambda})^{\ha^{+}} \end{array} \right.  =\label{2.29}\ee
$$ = \mbox{const}\, \psi_{n_{r}\beta_{r},
\ldots , n_{1}\beta_{1}} (m_{r} + \frac{ \langle n_{r}\beta_{r} , \mu
\rangle}{k}, \ldots ,m_{1} + \frac{ \langle n_{1}\beta_{1} , \mu
\rangle}{k} )\left |
\begin{array}[t]{l} \\
L_{\mu}(\hat{\Lambda})^{\ha^{+}} \end{array} \right.  ,$$
where $\mbox{const} \in \C^{\times}$ (cf. (\ref{2.104})). Since the
usual $\psi$-monomials are linear combinations of normalized
$\psi$-monomials, one can now conclude from (\ref{2.301}) and (\ref{2.21})
that
\be  L(\hat{\Lambda})^{\ha^{+}}= \rho(kQ)\cdot
\pi_{L(\hat{\Lambda})^{\ha^{+}}}\cdot
W(\hat{\Lambda}). \label{2.302}\ee
and
\be  L(\hat{\Lambda})^{\ha^{+}}_{kQ} =
\pi_{L(\hat{\Lambda})^{\ha^{+}}_{kQ}}\cdot
W(\hat{\Lambda})\label{2.30}\ee (cf. (\ref{1.6}) and (\ref{1.101})).
This close relationship with the principal subspace will be
instrumental in the subsequent  arguments.

\section{Quasiparticles}
\hspace{3ex}Recall that in the context of the principal subspace [GeI], we
restricted ourselves to vertex operators associated with the simple
roots because those are perfectly enough to generate the whole space
when acting on a highest weight vector $v(\hat{\Lambda})$ (Lemma 3.1
[GeI]). This property is inherited in the parafermionic picture. At
our convenience, we shall often be writing the $\psi$-monomials acting
on $L(\hat{\Lambda})^{\ha^{+}}$ as parafermionic counterparts (cf.
(\ref{2.104}) or (\ref{2.29})):\vspace{5mm}
\begin{Lemma}One has
\be L(\hat{\Lambda})^{\ha^{+}} = \mbox{\rm
Span}_{\C}\left\{ \rho(k\alpha)\cdot\pi_{L(\hat{\Lambda})^{\ha^{+}}}\cdot
x_{n_{r}\beta_{r}}(m_{r}) \cdots   \pi_{L(\hat{\Lambda})^{\ha^{+}}}\cdot
x_{n_{1}\beta_{1}}(m_{1})\cdot v(\hat{\Lambda})\right|\label{3.102}\ee
$$\left|x_{n_{r}\beta_{r}}(m_{r}) \cdots x_{n_{1}\beta_{1}}(m_{1})
\;\mbox{\rm a monomial from}\; U,\; \alpha \in Q \right\}.$$
Equivalently,
\be L(\hat{\Lambda})^{\ha^{+}}_{kQ} = \mbox{\rm
Span}_{\C}\left\{
\psi_{n_{r}\beta_{r}}(m_{r} + \frac{\langle n_{r}\beta_{r},
\Lambda + \sum_{p=1}^{r-1}n_{p}\beta_{p}\rangle}{k} ) \cdots
\psi_{n_{1}\beta_{1}}(m_{1} + \frac{\langle n_{1}\beta_{1},
\Lambda\rangle}{k} )\cdot v(\hat{\Lambda})\right|\label{3.1}\ee
$$\left|x_{n_{r}\beta_{r}}(m_{r}) \cdots x_{n_{1}\beta_{1}}(m_{1})
\;\mbox{\rm a monomial from}\; U  := U(\nb_{\alpha_{n}})\cdots
U(\nb_{\alpha_{1}})\right\}.$$
\label{Lem 3.1}\end{Lemma}

\noindent{Proof} Follows immediately from Lemma 3.1 [GeI], the
''surjectivity''  (\ref{2.302}) of the projection
$\pi_{L(\hat{\Lambda})^{\ha^{+}}}$ and the fact that normalized
monomials are linear combinations of usual monomials of the same
structure (cf. (\ref{2.104}) and (\ref{2.29})).\hfill
$\Box$\vspace{5mm}

   In other words, we are entitled again to dismiss all the nonsimple
roots and fix an order in the set of simple roots.\vspace{5mm}

\noindent{\bf Definition 3.1}  For every simple root $\alpha_{i},
\; 1\leq i
\leq n,$ and
positive integer $r,$ we shall say that the operator
$\psi_{r\alpha_{i}}(m)$ from (\ref{2.103}) represents a {\em
$\psi$-quasi-particle of color $i,$ charge $r$ and energy $-m
-\frac{r^{2}}{k}$.} If confusion with $x$-quasi-particles can not
arise, we shall skip the prefix $\psi$ and just talk about
quasi-particles; it will be clear from the context whether $\psi$ is
an operator on $L(\hat{\Lambda})^{\ha^{+}}$ or
$L(\hat{\Lambda})^{\ha^{+}}_{kQ}.$ Abusing language, we shall say
that $\psi_{\alpha_{i}}(m)$ is from
$\pi_{L(\hat{\Lambda})^{\ha^{+}}}\cdot U(\nb_{\alpha_{i}})$ and that
the $\psi$-monomials (normalized or not) considered in Lemma \ref{Lem
3.1} (\ref{3.102})
are from $\pi_{L(\hat{\Lambda})^{\ha^{+}}}\cdot U.$ Be aware that the
operator $\psi_{r\alpha_{i}}(m)$ is defined only on those $\mu$-weight
subspaces of $L(\hat{\Lambda})^{\ha^{+}}$ (resp.,
$L(\hat{\Lambda})^{\ha^{+}}),$ for which $m\in \Z +
\frac{\langle r\alpha_{i} , \mu \rangle}{k}.$\vspace{5mm}

   Starting from here,  we shall mostly work with ($\psi$-)  monomials
{}from  $\pi_{L(\hat{\Lambda})^{\ha^{+}}_{kQ}}\cdot U.$

   In complete analogy
with the $x$-monomials from $U$ ([GeI] Section 3), we shall say that a
$\psi$-monomial from $\pi_{L(\hat{\Lambda})^{\ha^{+}}}\cdot U$
\be \psi_{n_{r_{n}^{(1)},n}\alpha_{n}}(m_{r_{n}^{(1)},n})
\cdots \psi_{n_{1,n}\alpha_{n}}(m_{1,n}) \cdots \cdots
\psi_{n_{r_{1}^{(1)},1}\alpha_{1}}(m_{r_{1}^{(1)},1}) \cdots
\psi_{n_{1,1}\alpha_{1}}(m_{1,1})\label{3.3}\ee
and its normalized counterpart
\be \psi_{n_{r_{n}^{(1)},n}\alpha_{n} , \ldots
n_{1,n}\alpha_{n};\ldots \ldots
;n_{r_{1}^{(1)},1}\alpha_{1}, \ldots ,
n_{1,1}\alpha_{1}}(m_{r_{n}^{(1)},n},\ldots,
m_{1,n};\ldots \ldots ;m_{r_{1}^{(1)},1}, \ldots
,m_{1,1}),\label{3.4}\ee
are  of {\em color-charge-type}
\be (n_{r_{n}^{(1)},n},\ldots ,n_{1,n};\ldots ; n_{r_{1}^{(1)},1},
\ldots , n_{1,1}),\label{3.5}\ee
where
$$0 < n_{r^{(1)},i} \leq \cdots \leq n_{2,i} \leq n_{1,i} \leq
K,\;\sum_{p =1}^{r_{i}^{(1)}} n_{p,i} = r_{i},\; 1\leq i \leq n,$$
of {\em color-dual-charge-type}
\be (r_{n}^{(1)}, \ldots , r_{n}^{(K)}; \ldots ; r_{1}^{(1)},
\ldots , r_{1}^{(K)}),\label{3.6}\ee
where
$$ r_{i}^{(1)} \geq r_{i}^{(2)} \geq \cdots \geq r_{i}^{(K)} \geq 0,\; \sum_{t
=1}^{K}r_{i}^{(t)} = r_{i},\; K \in
\Z_{+},\;1\leq i \leq n,$$
and of {\em color-type} $(r_{n};
\ldots; r_{1}).$ We  shall also say that the corresponding generating
functions
\be \Psi_{n_{r_{n}^{(1)},n}\alpha_{n}}(z_{n_{r_{n}^{(1)},n}}) \cdots
\Psi_{n_{1,1}\alpha_{1}}(z_{1,1})\label{3.7}\ee
and
\be \Psi_{n_{r_{n}^{(1)},n}\alpha_{n},\ldots
,n_{1,1}\alpha_{1}}(z_{n_{r_{n}^{(1)},n}} ,\ldots , z_{1,1})\label{3.8}\ee
are of the above color-charge-type, color-dual-charge-type and
color-type.

   We would like to remind that no $\psi$-quasi-particles of charge
greater than $k-1$ will appear in our parafermionic space at level
$k$: the constraint (\ref{2.22}), (\ref{2.23}) asserts in particular
that for every color $i,\; 1\leq i \leq n,$ one has
\be  \Psi_{k\alpha_{i}}(z) = \mbox{const}\,\rho(k\alpha_{i}),\;
\mbox{const} \in \C^{\times} \label{3.13}\ee
as an operator on $L(\hat{\Lambda})^{\ha^{+}},$ i.e.,
\be  \Psi_{k\alpha_{i}}(z) = \mbox{const} \in \C^{\times} \label{3.14}\ee
as an operator  on the parafermionic space $L(\hat{\Lambda}
)^{\ha^{+}}_{kQ}.$ Note a curious
corollary  of this constraint: The
inverse of the identity (\ref{2.14}) implies that up to an invertible
operator, every
$x$-quasi-particle of charge $k$ and color $i$ acts on $L(\hat{\Lambda})$ as
an operator  from $U(\ha),$ i.e.,
\be x_{k\alpha_{i}}(m) = \rho(k\alpha_{i})\cdot h,\; h \in
U(\ha),\label{3.15}\ee
(recall that $ U(\ha)$ acts through the iterated coproduct
$\Delta^{k-1}).$

   Since our ambition is to exploit results from [GeI], it is really
more convenient to work with the $\ga \supset \ha$ subspace
$L(\hat{\Lambda})^{\ha^{+}}$ rather than with the parafermionic space
$L(\hat{\Lambda})^{\ha^{+}}_{kQ}$: the $\psi$-monomial
$\pi_{L(\hat{\Lambda})^{\ha^{+}}}\cdot x_{n_{r}\beta_{r}}(m_{r})
\cdots \pi_{L(\hat{\Lambda})^{\ha^{+}}}\cdot
x_{n_{1}\beta_{1}}(m_{1})$ and its normalized counterpart are of given
color-charge-type, color-dual-charge-type and color-type if their
$x$-counterpart $x_{n_{r}\beta_{r}}(m_{r})\cdots
x_{n_{1}\beta_{1}}(m_{1})$ is of these types. (Strictly speaking, the
$\psi$-monomials acting on $L(\hat{\Lambda} )^{\ha^{+}}_{kQ}$ do not
have $x$-counterparts because  $\pi_{L(\hat{\Lambda})^{\ha^{+}}_{kQ}}\cdot
x_{r\beta}(m)$ is not an operator on $L(\hat{\Lambda})^{\ha^{+}}_{kQ}.)$ Note
that given a $\psi$-monomial from
$\pi_{L(\hat{\Lambda})^{\ha^{+}}}\cdot U,$ we do not exactly know the
quasi-particle energies of its
$x$-counterpart unless we specify the weight subspace on which the
$\psi$-monomial acts  (cf.  (\ref{2.104}) and (\ref{2.29})).  If one
discusses only the type of a monomial, this is not necessary.

   In the same  flow of thoughts, it is clear how to
convey the  linear ordering ''$<$'' and the  partial ordering
''$\prec$'' from the set of $x$-monomials of a given color-type $(r_{n};
\ldots; r_{1})$ (cf. [GeI] Section 3) to  the set of $\psi$-monomials
(acting on $L(\hat{\Lambda})^{\ha^{+}}$) of a given  color-type: Set
\be \pi_{L(\hat{\Lambda})^{\ha^{+}}} \cdot
x_{n_{r}\beta_{r}}(m_{r})\cdots  \pi_{L(\hat{\Lambda})^{\ha^{+}}} \cdot
x_{n_{1}\beta_{1}}(m_{1}) <
\pi_{L(\hat{\Lambda})^{\ha^{+}}} \cdot
x_{n'_{r}\beta_{r}}(m'_{r})\cdots \pi_{L(\hat{\Lambda})^{\ha^{+}}}
\cdot x_{n'_{1}\beta_{1}}(m'_{1})\label{3.9}\ee
and
\be \pi_{L(\hat{\Lambda})^{\ha^{+}}} \cdot x_{n_{r}\beta_{r}}(m_{r})\cdots
x_{n_{1}\beta_{1}}(m_{1}) < \pi_{L(\hat{\Lambda})^{\ha^{+}}} \cdot
x_{n'_{r}\beta_{r}}(m'_{r})\cdots
x_{n'_{1}\beta_{1}}(m'_{1})\label{3.10}\ee
if
$$x_{n_{r}\beta_{r}}(m_{r})\cdots
x_{n_{1}\beta_{1}}(m_{1})<x_{n'_{r}\beta_{r}}(m'_{r})\cdots
x_{n'_{1}\beta_{1}}(m'_{1}).$$
Define analogously ''$\prec$'' for $\psi$-monomials of a given
color-type. These definitions can obviously be rewritten, replacing
the projections with
the corresponding $\psi$-operators according to (\ref{2.104}) and
(\ref{2.29}). Moreover, we do not have to
specify the
weight subspace on which  the $\psi$-monomials  act: Recall from [GeI]
Section 3 that  the quassi-particle energies
affect the ordering only if the color-charge-types (and for that
matter, the color-dual-charge-types) are the same, in which case the
shift of the indices of the corresponding $x$-monomials will be the
same for the two compared $\psi$-monomials. Keep in mind  that ''$\prec$''
implies  ''$<$''   but not vice versa.

   We are now in position to  explain why  working
with usual or  with normalized $\psi$-monomials is essentially the
same. Roughly speaking, the transformation matrix between them is ''upper
triangular''. More precisely, from the definitions (\ref{2.24}),
(\ref{2.103}),(\ref{2.27}) and the identities  (\ref{2.104}),
(\ref{2.29}), we obtain for a given (usual) monomial
{}from $\pi_{L(\hat{\Lambda})^{\ha^{+}}} \cdot U$ that
\be  \pi_{L(\hat{\Lambda})^{\ha^{+}}} \cdot
x_{n_{r}\beta_{r}}(m_{r})\cdots  \pi_{L(\hat{\Lambda})^{\ha^{+}}} \cdot
x_{n_{1}\beta_{1}}(m_{1}) = \label{3.11}\ee
$$=\mbox{const}\; \pi_{L(\hat{\Lambda})^{\ha^{+}}} \cdot
x_{n_{r}\beta_{r}}(m_{r})\cdots x_{n_{1}\beta_{1}}(m_{1})
+ \begin{array}[t]{l} \mbox{a linear combination of normalized}\\ \mbox{
monomials of the same color-charge-}\\\mbox{-type and greater in the
ordering ''$\prec$'',}\end{array}$$
where $\mbox{const} \in \C^{\times}.$
Conversely, for any given normalized monomial from
$\pi_{L(\hat{\Lambda})^{\ha^{+}}} \cdot U,$ we have
\be  \pi_{L(\hat{\Lambda})^{\ha^{+}}} \cdot
x_{n_{r}\beta_{r}}(m_{r})\cdots x_{n_{1}\beta_{1}}(m_{1}) =
\label{3.12}\ee
$$= \mbox{const}\;  \pi_{L(\hat{\Lambda})^{\ha^{+}}} \cdot
x_{n_{r}\beta_{r}}(m_{r})\cdots  \pi_{L(\hat{\Lambda})^{\ha^{+}}} \cdot
x_{n_{1}\beta_{1}}(m_{1})+ \begin{array}[t]{l} \mbox{a linear
combination of (usual)}\\  \mbox{
monomials of the same color-}\\\mbox{-charge-type and
greater}\\\mbox{in the ordering ''$\prec$'',}\end{array}$$
$\mbox{const} \in \C^{\times}.$

\section{Quasi-particle basis for the parafermionic space}
\hspace{3ex}Recall that we are working with level $k \in  \Z_{+},\; k
\geq 2,$ and a highest weight
\be \hat{\Lambda} :=
k_{0}\hat{\Lambda}_{0} + k_{j}\hat{\Lambda}_{j} = k\hat{\Lambda}_{0} +
\Lambda,\;\;\mbox{where}\; \Lambda :=  k_{j}\Lambda_{j},\label{4.1}\ee
for some $j,\; 1 \leq j \leq n;$ $k_{0}, k_{j}
\in \N $ and $k_{0} + k_{j} = k$ (cf. Section 1).  Our level $k$
spaces are all realized in the tensor product of level one modules
given by the homogeneous vertex operator construction ([GeI] Section
2). In order to treat the level one ingredients on equal footing,
we introduced
\be  j_{t} :=  \left\{ \begin{array}{l} \ 0\;\; \mbox{for} \;
0 < t \leq  k_{0} \\  j\;\;  \mbox{for} \;
k_{0}< t \leq k = k_{0}+ k_{j}\end{array}
\right.\label{4.2}\ee
and then the highest weight vector of all the modules under
consideration $L(\hat{\Lambda}),$ $ W(\hat{\Lambda}),$
$L(\hat{\Lambda})^{\ha^{+}}$ and
$L(\hat{\Lambda})^{\ha^{+}}_{kQ}$ was defined as
\be
v(\hat{\Lambda}) := v(\hat{\Lambda}_{j_{k}})\otimes \cdots \otimes
v(\hat{\Lambda}_{j_{1}}) = \label{4.3}\ee
$$=\underbrace{v(\hat{\Lambda}_{j})\otimes\cdots \otimes
v(\hat{\Lambda}_{j})}_{k_{j}\;factors} \otimes
\underbrace{v(\hat{\Lambda}_{0})\otimes\cdots \otimes
v(\hat{\Lambda}_{0})}_{k_{0}\;factors},$$
where $v(\hat{\Lambda}_{j_{t}})$ is the highest weight vector of the
level one module in the $t^{th}$ tensor slot (counted from right to
left).

   The insightful reader has probably guessed already what is our
basis-candidate for the level $k$ parafermionic space
$L(\hat{\Lambda})^{\ha^{+}}_{kQ}$ -- the most intuitive choice is of
course the $\pi_{L(\hat{\Lambda})^{\ha^{+}}_{kQ}}$-projection of this
particular subset of the basis for $W(\hat{\Lambda})$ (from [GeI]
Section 5), which contains only vectors generated by monomials with no
quasi-particles of charge $k.$ We shall rewrite for completeness the
full definition of the prototype -- the basis-generating set
$\frak{B}_{W(\hat{\Lambda})}$ from Definition 5.1 [GeI] -- with the
charge of the quasi-particles bounded by $k-1$ rather than $k,$ and
call it $\frak{B}_{W(\hat{\Lambda})}^{(k-1)}.$ We shall not comment
here on the origin and naturality of this incomprehensible at first
sight affluence of parameters.  The frustrated readers are referred to
the Introduction of [GeI] for simple particular cases and to the
Introduction and
Section 5 here for much easier to grasp character formulas associated
with these bases. We only emphasize that the
basis-generating set of monomials is a disjoint union along
color-charge-types (\ref{3.5}) or, equivalently, along the
corresponding color-dual-charge-types (\ref{3.6}), with upper bound
for the charges $K := k-1$ (as opposed to $K = k$ in the principal
subspace picture):\newpage

\noindent {\bf Definition 4.1} Fix a highest weight $\hat{\Lambda}$ as in
(\ref{4.1}). Set

\be \frak{B}_{W(\hat{\Lambda})}^{(k-1)} :=
\bigsqcup_{\begin{array}{c}{\scriptstyle 0 \leq
n_{r_{n}^{(1)},n} \leq \cdots \leq  n_{1,n}\leq k-1}\\{\scriptstyle
\cdots\cdots\cdots}\\{\scriptstyle 0 \leq
n_{r_{1}^{(1)},1} \leq \cdots \leq  n_{1,1}\leq
k-1}\end{array}}\left(\mbox{or,
equivalently,}\;\;\bigsqcup_{\begin{array}{c}{\scriptstyle
r_{n}^{(1)}\geq  \cdots \geq  r_{n}^{(k-1)}\geq 0}\\{\scriptstyle
\cdots\cdots\cdots}\\{\scriptstyle
r_{1}^{(1)} \geq \cdots \geq  r_{1}^{(k-1)}\geq
0}\end{array}}\right)
\label{4.4}\ee
$$\left\{
x_{n_{r_{n}^{(1)},n}\alpha_{n}}(m_{r_{n}^{(1)},n})
\cdots x_{n_{1,n}\alpha_{n}}(m_{1,n}) \cdots \cdots
x_{n_{r_{1}^{(1)},1}\alpha_{1}}(m_{r_{1}^{(1)},1}) \cdots
x_{n_{1,1}\alpha_{1}}(m_{1,1})\begin{array}{c}\\ \\ \\
\end{array}\right|$$
$$\left|\begin{array}{ll} m_{p,i} \in \Z,\; 1\leq i \leq n, \; 1 \leq
p \leq r_{i}^{(1)};&\\m_{p,i} \leq
\sum_{q=1}^{r_{i-1}^{(1)}} \mbox{min}\,\{n_{p,i}, n_{q,i-1}\}  -
\sum_{t=1}^{n_{p,i}}\delta_{i,j_{t}}
-\sum_{p>p'>0} 2\mbox{min}\, \{n_{p,i}, n_{p',i}\} -
n_{p,i};&\\m_{p+1,i} \leq m_{p,i}
-2n_{p,i}\;\;\mbox{for}\;\; n_{p+1,i} = n_{p,i}& \end{array}\right\}, $$
where $r_{0}^{(1)} := 0$ and  $j_{t}$ was introduced  in (\ref{4.2}).
Define
\be  \frak{B}_{L(\hat{\Lambda})^{\ha^{+}}} :=
\left\{\pi_{L(\hat{\Lambda})^{\ha^{+}}} \cdot
x_{n_{r}\beta_{r}}(m_{r})\cdots  \pi_{L(\hat{\Lambda})^{\ha^{+}}} \cdot
x_{n_{1}\beta_{1}}(m_{1})\right| \label{4.6}\ee
$$\left| x_{n_{r}\beta_{r}}(m_{r})\cdots
x_{n_{1}\beta_{1}}(m_{1}) \in
\frak{B}_{W(\hat{\Lambda})}^{(k-1)}\right\}$$
and
\be  \frak{B}_{L(\hat{\Lambda})^{\ha^{+}}_{kQ}} :=
\left\{\psi_{n_{r}\beta_{r}}(m_{r} + \frac{\langle n_{r}\beta_{r},
\Lambda + \sum_{p=1}^{r-1}n_{p}\beta_{p}\rangle}{k} ) \cdots
\psi_{n_{1}\beta_{1}}(m_{1} + \frac{\langle n_{1}\beta_{1},
\Lambda\rangle}{k} )    \right| \label{4.5}\ee
$$\left| x_{n_{r}\beta_{r}}(m_{r})\cdots
x_{n_{1}\beta_{1}}(m_{1}) \in
\frak{B}_{W(\hat{\Lambda})}^{(k-1)}\right\}$$
(cf. (\ref{1.6}), (\ref{1.8}) and (\ref{1.101})).  \label{Def
4.1}\vspace{5mm}

\noindent{\it Example 4.1} Consider $\g = sl(3),$ i.e., $n=2$ and the vacuum
highest weight $\hat{\Lambda} = 2\hat{\Lambda}_{0}$ at level $k=2.$
Denote for brevity the monomial
$$\pi_{L(2\hat{\Lambda}_{0})^{\ha^{+}}}
\cdot x_{\alpha_{2}}(s) \cdots \pi_{L(2\hat{\Lambda}_{0})^{\ha^{+}}}
\cdot
x_{\alpha_{1}}(t)$$ by $(s_{\alpha_{2}} \ldots t_{\alpha_{1}}).$ For
the first few energy levels (the eigenvalues of the scaling operator
$D- D^{\ha}$ under the adjoint action), Table 1 in  the Appendix lists the
elements of $\frak{B}_{L(2\hat{\Lambda}_{0})^{\ha^{+}}}$ of
color-types $(1;2)$
and $(2;2).$ This Table is of course a copy of Table 1 [GeI],
Appendix (listing the corresponding elements of
$\frak{B}_{W(2\hat{\Lambda}_{0})}^{(1)} \cong
\frak{B}_{W(\hat{\Lambda}_{0})})$  with the entries in the column
''energy'' shifted  according to (\ref{2.105}) and (\ref{2.10}).
\vspace{5mm}

\noindent {\it Example 4.2 } Let again $\g = sl(3)$ but consider the vacuum
highest weight $\hat{\Lambda} = 3\hat{\Lambda}_{0}$  at level
$k=3.$ Similarly to the previous example, denote the quasi-particle
monomial $$\pi_{L(3\hat{\Lambda}_{0})^{\ha^{+}}}
\cdot x_{s'\alpha_{2}}(s) \cdots \pi_{L(3\hat{\Lambda}_{0})^{\ha^{+}}}
\cdot
x_{t'\alpha_{1}}(t)$$ by $(s_{s'\alpha_{2}} \ldots t_{t'\alpha_{1}}).$
For the first few energy levels, Table 2 in  the
Appendix lists the elements of
$\frak{B}_{L(3\hat{\Lambda}_{0})^{\ha^{+}}}$ of color-types
$(1;2)$ and $(2;2).$ This Table is a copy of Table 2 [GeI],
Appendix (listing the corresponding elements of
$\frak{B}_{W(3\hat{\Lambda}_{0})}^{(2)} \cong
\frak{B}_{W(2\hat{\Lambda}_{0})})$  with the entries in the column
''energy'' shifted  according to (\ref{2.106}) and (\ref{2.104}).
\vspace{5mm}

   We  begin as usual with a proof of the spanning property of our
basis-candidate.\vspace{5mm}

\begin{Theorem} Let $\hat{\Lambda}$ be  a highest weight as in (\ref{4.1}) and
$v(\hat{\Lambda})$ be the corresponding highest weight vector (\ref{4.3}).
Then the set  $\left\{\rho(k\alpha) \cdot b
\cdot v(\hat{\Lambda})\,\left| b \in
\frak{B}_{L(\hat{\Lambda})^{\ha^{+}}},\; \alpha \in Q \right\}
\right.$ spans the $\ga \supset \ha$ subspace $L(\hat{\Lambda})^{\ha^{+}}.$
Equivalently, the set $
\left\{b \cdot v(\hat{\Lambda})\,\left| b \in
\frak{B}_{L(\hat{\Lambda})^{\ha^{+}}_{kQ}}\right\} \right.$ spans the
parafermionic space $L(\hat{\Lambda})^{\ha^{+}}_{kQ}.$ \label{The 4.1}
\end{Theorem}

\noindent{\em Proof} We shall prove the statement for
$L(\hat{\Lambda})^{\ha^{+}}$ (the statement for
$L(\hat{\Lambda})^{\ha^{+}}_{kQ}$ follows immediately from
(\ref{2.104}) and the definition (\ref{1.9})).  In view of Lemma
\ref{Lem 3.1}, it suffices to show that every vector $b\cdot
v(\hat{\Lambda}),$ $b$ an usual $\psi$-monomial from
$\pi_{L(\hat{\Lambda})^{\ha^{+}}} \cdot U,$ is a linear combination of
vectors from the proposed set.

   Suppose that $b \not\in \rho(kQ)\cdot\frak{B}_{L(\hat{\Lambda})^{\ha^{+}}}.$
Due to the constraint (\ref{3.13}), the
(nonzero) quasi-particles of charge $k$ in $b$ (if any!) commute with
all the other quasi-particles  and can be moved (``exiled'') all the way to the
left. Using the ''upper triangular'' transformation (\ref{3.11}), normalize
the remaining $\psi$-monomial (which is by assumption $\not\in
\frak{B}_{L(\hat{\Lambda})^{\ha^{+}}})$ and
apply Remark 5.1 [GeI] -- the strong form of the spanning Theorem 5.1
[GeI] --  for
the obtained
$x$-monomials (new quasi-particles of charge $k$ might be generated in
the process!). Then  switch back from normalized to usual
$\psi$-monomials using the inverse transformation (\ref{3.12}) and
return the ''exiled''
quasi-particles of charge $k$ to their old places. The
newly obtained (usual) $\psi$-monomials  from
$\pi_{L(\hat{\Lambda})^{\ha^{+}}} \cdot U$ have the same color-type and
index-sum as $b$ and
moreover, since  $b \not\in
\rho(kQ)\cdot\frak{B}_{L(\hat{\Lambda})^{\ha^{+}}},$ they are
all greater than $b$ in the ordering $''\prec$'' (by Remark 5.1
[GeI]). Since there are only finitely many such $\psi$-monomials which
do not annihilate $v(\hat{\Lambda}),$ the statement follows by
induction.\hfill $\Box$ \vspace{5mm}

\noindent {\it Remark 4.1} The proof implies that the exact analog of Remark
5.1 [GeI] in the current setting is also true: Every vector $b\cdot
v(\hat{\Lambda}),$ $b$ a  usual $\psi$-monomial from
$\pi_{L(\hat{\Lambda})^{\ha^{+}}}
\cdot U,$ $ b \not\in \rho(kQ)\cdot\frak{B}_{L(\hat{\Lambda})^{\ha^{+}}},$
is a linear
combination of vectors of the form $b'\cdot
v(\hat{\Lambda}),$ $ b' \in
\rho(kQ)\cdot\frak{B}_{L(\hat{\Lambda})^{\ha^{+}}},$
$b' \succ b,$
with $b'$ and $b$ having the same color-type and total
index-sum. \vspace{5mm}

   We proceed with the independence. For the familiar  highest weight
$\hat{\Lambda} =
\sum_{t=1}^{k} \hat{\Lambda}_{j_{t}}$ (cf. (\ref{4.1}) and
(\ref{4.2})), define
\be \hat{\dot{\Lambda}}: =
\sum_{t=1}^{k-1} \hat{\Lambda}_{j_{t}} = \hat{\Lambda} -
\hat{\Lambda}_{j_{k}},\label{4.9}\ee
i.e., $\hat{\dot{\Lambda}} = (k-1)\hat{\Lambda}_{0} + \dot{\Lambda},$
where $\dot{\Lambda} = \Lambda - \Lambda_{j_{k}}$ or
equivalently,  $\dot{\Lambda}  = (k_{j} -1)
\Lambda_{j} \in P$ if $k_{j} >
0$ and $\dot{\Lambda} = 0$ otherwise.
Note that $\hat{\dot{\Lambda}}$ is of the same type (\ref{4.1}) as
$\hat{\Lambda}$ (i.e., the results in [GeI] hold for
$\hat{\dot{\Lambda}}$) and moreover, by the very definition
(\ref{4.3}), one has for the corresponding highest weight vectors
$v(\hat{\Lambda}) = v(\hat{\Lambda}_{j_{k}}) \otimes
v(\hat{\dot{\Lambda}}).$ Summoning the projection
$\pi_{U(\ha^{-})\cdot v(\hat{\Lambda}_{j})}$ from [GeI] (2.12), one
easily checks that
\be  \left(\pi_{U(\ha^{-})\cdot v(\hat{\Lambda}_{j_{k}})}\otimes
\underbrace{\mbox{id}\otimes\cdots \otimes
\mbox{id}}_{k-1\,factors}\right) \cdot \frak{B}_{W(\hat{\Lambda})}\cdot
v(\hat{\Lambda}) =\label{4.10}\ee
$$=  \left(\pi_{U(\ha^{-})\cdot v(\hat{\Lambda}_{j_{k}})}\otimes
\underbrace{\mbox{id}\otimes\cdots \otimes
\mbox{id}}_{k-1\,factors}\right) \cdot
\frak{B}_{W(\hat{\Lambda})}^{(k-1)}\cdot v(\hat{\Lambda})  =
v(\hat{\Lambda}_{j_{k}}) \otimes
\frak{B}_{W(\hat{\dot{\Lambda}})}\cdot v(\hat{\dot{\Lambda}}),$$ (cf.
Definition 4.1 and [GeI] Definition 5.1). This is
because any nontrivial $x$-quasi-particle action on the leftmost
tensor factor vanishes under the above projection due to the term
$e^{\beta}$ in the quasi-particle (cf. (\ref{2.1}), (\ref{2.101})).
But $x$-quasi-particles of charge $k$ are zero unless they act on all
the $k$ tensor factors.

   In the independence argument below, we shall employ the
independence of the vectors from the set
$\frak{B}_{W(\hat{\dot{\Lambda}})}\cdot v(\hat{\dot{\Lambda}})$
(rather than $\frak{B}_{W(\hat{\Lambda})}\cdot v(\hat{\Lambda})!$)
which was proven in Theorem 5.2 [GeI].
\vspace{5mm}

\begin{Theorem} Let $\hat{\Lambda}$ be  again a highest weight as in
(\ref{4.1}) and $v(\hat{\Lambda})$ be the corresponding highest weight
vector (\ref{4.3}).  Then the set $
\left\{\rho(k\alpha)\cdot b
\cdot v(\hat{\Lambda})\,\left| b \in
\frak{B}_{L(\hat{\Lambda})^{\ha^{+}}},\; \alpha \in Q\right\} \right.$
is indeed a basis for the $\ga \supset \ha$ subspace
$L(\hat{\Lambda})^{\ha^{+}}.$ Equivalently, the set $
\left\{ b
\cdot v(\hat{\Lambda})\,\left| b \in
\frak{B}_{L(\hat{\Lambda})^{\ha^{+}}_{kQ}} \right\} \right.$  is
a basis for the parafermionic space $L(\hat{\Lambda})^{\ha^{+}}_{kQ}.$
\label{The 4.2}\end{Theorem}

\noindent{\em Proof} It suffices to  prove the independence of the vectors
in the set $\rho(kQ)\cdot\frak{B}_{L(\hat{\Lambda})^{\ha^{+}}} \cdot
v(\hat{\Lambda}).$ Let us first show that it would follow from the
independence of the vectors in  the set $v(\hat{\Lambda}_{j_{k}})
\otimes
\frak{B}_{W(\hat{\dot{\Lambda}})}\cdot v(\hat{\dot{\Lambda}})$ {\em
over the ring} $U(\ha^{-})$ (recall that $U(\ha)$ acts through the
$(k-1)$-iterate $\Delta^{k-1}$ of the standard coproduct $\Delta).$

   Suppose that there is a (nontrivial irreducible) linear relation
among vectors from  $\rho(kQ)\cdot\frak{B}_{L(\hat{\Lambda})^{\ha^{+}}}
\cdot v(\hat{\Lambda}).$ Without loss of generality, one can assume
that the relation does not contain factors from $\rho(kQ_{-})$ (recall
{}from [GeI] Preliminaries that $Q_{-} \subset Q$ is the semigroup
generated by the negatives of the simple roots) and that at least one
vector is from $\frak{B}_{L(\hat{\Lambda})^{\ha^{+}}}\cdot
v(\hat{\Lambda}).$ An induction on the number of monomials involved
shows that this can be easily achieved by multiplying the relation
with appropriate invertible operators from $\rho(kQ).$

   Normalize the
$\psi$-monomials using the ''upper triangular'' relation
(\ref{3.11}): Observe that by (\ref{2.25}), a vector
$\pi_{L(\hat{\Lambda})^{\ha^{+}}} \cdot b \cdot v(\hat{\Lambda}),$ $b$
a monomial
{}from $U,$ equals up to a nonzero constant $b \cdot v(\hat{\Lambda}),$ plus
a linear combination of vectors of the form
$h'\cdot b'\cdot v(\hat{\Lambda}),$ where $h' \in U(\ha^{-}),$ $b'$ is a
monomial from $U$ of the same color-type but of less index-sum than
$b,$ and in addition, $b' \succ b.$ Implement  the strong
form of the spanning Theorem 5.1 [GeI] (cf. Remark 5.1 [GeI]) for the
obtained $x$-monomials and apply the projection
\be \pi_{U(\ha^{-})\cdot
v(\hat{\Lambda}_{j_{k}})}\otimes \underbrace{\mbox{id}\otimes\cdots \otimes
\mbox{id}}_{k-1\,factors}\label{4.11}\ee
to the relation. (Note that it  annihilates  all the
vectors containing a
factor $\rho (k\alpha),\; \alpha \in Q.)$  Due to (\ref{4.10}), the result
is a linear relation among vectors from
$v(\hat{\Lambda}_{j_{k}}) \otimes \frak{B}_{W(\hat{\dot{\Lambda}})}\cdot
v(\hat{\dot{\Lambda}}),$ with coefficients in $U(\ha^{-}).$ The  relation is
{\em nontrivial}: Among all the vectors from
$\frak{B}_{L(\hat{\Lambda})^{\ha^{+}}} \cdot v(\hat{\Lambda})$ in our
initial relation (we have seen that without loss of generality, there
is always at least one such vector),
let $b\cdot v(\hat{\Lambda})$ be the one whose $\psi$-monomial $b$ is
smallest in the  linear
ordering ''$<$''.  Then the projection (\ref{4.11}) of the
$x$-counterpart of $b$ acting on $
v(\hat{\Lambda}),$ is nonzero (by (\ref{4.10}) and Theorem 5.2 [GeI])
and is
present in the final relation, because ''$\prec$'' implies ''$<$''.

   In order to complete the proof, it remains to show that  a linear
relation among vectors from
$v(\hat{\Lambda}_{j_{k}}) \otimes \frak{B}_{W(\hat{\dot{\Lambda}})}\cdot
v(\hat{\dot{\Lambda}}),$ with coefficients in $U(\ha^{-})$ is
impossible. We shall reach a contradiction  with the
independence of the vectors from $\frak{B}_{W(\hat{\dot{\Lambda}})}\cdot
v(\hat{\dot{\Lambda}})$ by restricting the relation to an appropriate
homogeneous subspace where
the action of the polynomial algebra $U(\ha^{-})$ (given by
$\Delta^{k-1})$ is ''squeezed''  to the leftmost tensor slot: Suppose
that $m \in \N$ is the
maximal degree of the polynomials from $U(\ha^{-})$ which are
coefficients in our relation. Apply the projection
\be \pi_{U^{m}(\ha^{-})\cdot
v(\hat{\Lambda}_{j_{k}})}\otimes \underbrace{\mbox{id}\otimes\cdots \otimes
\mbox{id}}_{k-1\,factors},\label{4.12}\ee
where $\pi_{U^{m}(\ha^{-})\cdot v(\hat{\Lambda}_{j})}$ was introduced
in (2.12) [GeI]. Since the vectors from the set
$v(\hat{\Lambda}_{j_{k}}) \otimes
\frak{B}_{W(\hat{\dot{\Lambda}})}\cdot v(\hat{\dot{\Lambda}})$ are
nonzero (Theorem 5.2 [GeI]), the number $m$ is maximal with the
property that the projection (\ref{4.12}) does not annihilate all the
vectors in the relation. The reason we would like to apply the
projection (\ref{4.12}) to our relation is quite transparent: It is
easily seen from the explicit form of the iterated coproduct
$\Delta^{k-1}$ that we shall thus obtain a {\it nontrivial} relation
among vectors of the form $h\cdot v(\hat{\Lambda}_{j_{k}}) \otimes
b\cdot v(\hat{\dot{\Lambda}}),$ where $h \in U^{m}(\ha^{-})$ and $b
\in \frak{B}_{W(\hat{\dot{\Lambda}})}.$ But this contradicts the
independence of the vectors from the set
$\frak{B}_{W(\hat{\dot{\Lambda}})}\cdot v(\hat{\dot{\Lambda}})$ which
was proven in Theorem 5.2 [GeI].\hfill $\Box$ \vspace{5mm}

\section{Character formulas}
\hspace{3ex}We can finally reward  ourselves and enjoy   the
entertaining world of characters associated with the above bases.

   Consider an arbitrary monomial
\be \pi_{L(\hat{\Lambda})^{\ha^{+}}} \cdot
x_{n_{r}\beta_{r}}(m_{r})\cdots \pi_{L(\hat{\Lambda})^{\ha^{+}}} \cdot
x_{n_{1}\beta_{1}}(m_{1})\label{5.1}\ee
{}from $\pi_{L(\hat{\Lambda})^{\ha^{+}}}
\cdot U$ of color-dual-charge-type
\be (r_{n}^{(1)}, \ldots , r_{n}^{(k-1)}; \ldots ; r_{1}^{(1)},
\ldots , r_{1}^{(k-1)}),\label{5.2}\ee
where $$ r_{i}^{(1)} \geq r_{i}^{(2)} \geq \cdots \geq r_{i}^{(k-1)}
\geq 0,\; \sum_{t =1}^{k-1}r_{i}^{(t)} = r_{i},\;1\leq i \leq n,$$
and hence, of color-type $(r_{n};\ldots; r_{1})$ (cf. (\ref{3.6})).
Set $p_{i}^{(s)}$ to be the number of quasi-particles of color $i$ and
charge $s$ in our monomial, i.e., $p_{i}^{(s)} := r_{i}^{(s)} -
r_{i}^{(s+ 1)},$ $1 \leq s \leq k-2$ and $p_{i}^{(k - 1)} := r_{i}^{(k
- 1)},$ hence $r_{i} = \sum_{s =1}^{k-1}s p_{i}^{(s)},$ $i = 1, \ldots,
n.$ Then according to (\ref{2.28}), (\ref{2.29}) and (\ref{3.11}), the
$D - D^{\ha}$-eigenvalue of the vector
\be \pi_{L(\hat{\Lambda})^{\ha^{+}}} \cdot
x_{n_{r}\beta_{r}}(m_{r})\cdots \pi_{L(\hat{\Lambda})^{\ha^{+}}} \cdot
x_{n_{1}\beta_{1}}(m_{1})\cdot v(\hat{\Lambda})\label{5.3}\ee
equals
\be -\sum_{l= 1}^{r} m_{l}  - \frac{1}{2k}\langle
\sum_{i=1}^{n}(\sum_{s=1}^{k-1}s p_{i}^{(s)})\alpha_{i} ,
\sum_{i=1}^{n}(\sum_{s=1}^{k-1}s p_{i}^{(s)}) \alpha_{i}\rangle  -
\frac{1}{k}\langle
\sum_{i=1}^{n}(\sum_{s=1}^{k-1}s p_{i}^{(s)}) \alpha_{i}
,\Lambda \rangle = \label{5.4}\ee
$$ = -\sum_{l= 1}^{r} m_{l}  - \frac{1}{2k}\langle
\sum_{i=1}^{n}r_{i}\alpha_{i} ,
\sum_{i=1}^{n}r_{i}\alpha_{i}\rangle  -  \frac{1}{k}\langle
\sum_{i=1}^{n}r_{i}\alpha_{i}
,\Lambda \rangle  = $$
$$= -\sum_{l= 1}^{r} m_{l} -
\frac{1}{2}\sum_{\begin{array}{c} {\scriptstyle l,m= 1, \ldots,
n}\\{\scriptstyle s,t = 1, \ldots , k - 1}\end{array}}\; A_{lm}C^{st}
p_{l}^{(s)}p_{m}^{(t)} - \frac{k_{j}}{k} \sum _{s=1}^{k-1}s
p_{j}^{(s)} = $$ $$= -\sum_{l= 1}^{r} m_{l} - \frac{1}{k} (r_{1}^{2} +
\cdots + r_{n}^{2} - r_{1}r_{2} - \cdots - r_{n-1}r_{n}) -
\frac{k_{j}}{k} r_{j},$$ where $(A_{lm})_{l,m =1}^{n}$ is the Cartan
matrix of $\g$ and $C^{st}:= \frac{st}{k},\; 1\leq s,t \leq k-1.$ In
view of Definition 4.1, (\ref{4.9}) and Theorem
\ref{The 4.2}, this formula provides the correction term needed to
obtain $\mbox{Tr}\,q^{D- D^{\ha}}\left | \begin{array}[t]{l} \\
L(\hat{\Lambda})^{\ha^{+}}_{kQ} \end{array} \right.$ from
$\mbox{Tr}\,q^{D}\left | \begin{array}[t]{l} \\
\frak{B}_{W(\hat{\Lambda})}^{(k-1)}\cdot
v(\hat{\Lambda})\end{array} \right. = \mbox{Tr}\,q^{D}\left |
\begin{array}[t]{l} \\ W(\hat{\dot{\Lambda}})\end{array} \right.$ (cf.
also (\ref{4.10})). Recall from [GeI] (5.27) that
\be \mbox{\rm Tr}\,q^{D} \left| \begin{array}[t]{l} \\
W(\hat{\dot{\Lambda}}) \end{array} \right.
=\label{5.5}\ee
$$= \sum_{\begin{array}{c}{\scriptstyle
p_{1}^{(1)},\ldots , p_{1}^{(k-1)} \geq 0}\\{\scriptstyle \ldots \ldots
\ldots\ldots}\\{\scriptstyle p_{n}^{(1)},\ldots
, p_{n}^{(k-1)} \geq 0 }\end{array}}
\frac{q^{\frac{1}{2}\sum_{l,m =
1, \ldots, n}^{s,t = 1,\ldots , k-1}\;A_{lm} B^{st}
p_{l}^{(s)}p_{m}^{(t)}}}{\prod_{i=1}^{n} \prod_{s=1}^{k-1}
(q)_{p_{i}^{(s)}} }\;q^{\sum_{s=k_{0}+1}^{k-1} (s-k_{0})
p_{j}^{(s)}}$$ where $B^{st} := \mbox{min}\{s,t\},$ $1\leq s,t
\leq k-1.$ But one can immediately check that
\be B^{st} - C^{st} = A_{st}^{(-1)},\; 1\leq s,t \leq
k-1,\label{5.6}\ee where $(A_{st}^{(-1)})_{s,t=1}^{k-1}$ is the
inverse of the Cartan matrix of $sl(k,\C).$ The last three identities
and Theorem \ref{The 4.2} therefore imply
\be  \mbox{Tr}\,q^{D- D^{\ha}}\left | \begin{array}[t]{l} \\
L(\hat{\Lambda})^{\ha^{+}}_{kQ}
\end{array} \right. = \mbox{Tr}\,q^{D- D^{\ha}}\left | \begin{array}[t]{l} \\
L(k_{0}\hat{\Lambda}_{0} + k_{j}\hat{\Lambda}_{j})^{\ha^{+}}_{kQ}
\end{array} \right. = \label{5.7}\ee
$$= \sum_{\begin{array}{c}{\scriptstyle
p_{1}^{(1)},\ldots , p_{1}^{(k-1)} \geq 0}\\{\scriptstyle \ldots \ldots
\ldots\ldots}\\{\scriptstyle p_{n}^{(1)},\ldots
, p_{n}^{(k-1)} \geq 0 }\end{array}}
\frac{q^{\frac{1}{2}\sum_{l,m =
1, \ldots, n}^{s,t = 1,\ldots , k-1}\;A_{lm} A_{st}^{(-1)}
p_{l}^{(s)}p_{m}^{(t)}}}{\prod_{i=1}^{n} \prod_{s=1}^{k-1}
(q)_{p_{i}^{(s)}} }\;q^{\sum_{s=k_{0}+1}^{k-1} (s-k_{0})
p_{j}^{(s)}- \frac{k_{j}}{k} \sum _{s=1}^{k-1}s p_{j}^{(s)}}.$$
If we restrict our attention to  the vacuum module ($\hat{\Lambda} = k
\hat{\Lambda}_{0},$ that is $k_{0} = k$ and $k_{j} = 0),$ this formula
is the $sl(n+1,\C)$-case of the Kuniba-Nakanishi-Suzuki conjecture
[KNS] (a dilogarithm proof of this particular case was announced in
[Kir]).

   Note that for  a given weight $\mu \in P,$  such that $\mu -
\Lambda \in Q,$ one obtains the $(D - D^{\ha})$-character of
the weight subspace
$L_{\mu}(\Lambda)^{\ha^{+}}_{kQ} \cong L_{\mu}(\Lambda)^{\ha^{+}}$  if
the additional restriction
\be \mu =   \Lambda +  \sum_{i=1}^{n} r_{i}\alpha_{i} \;\;{\rm
mod}\;\;kQ  = \Lambda +
\sum_{i=1}^{n}(\sum_{s=1}^{k-1}s p_{i}^{(s)})\alpha_{i} \;\;{\rm
mod}\;\;kQ\label{5.8}\ee is imposed in the above formula (\ref{5.7}).
As far as the standard $q$-character is concerned, we know from
(\ref{2.201}) and (\ref{1.4}) that for $\mu \in P,$
\be  \mbox{\rm Tr}\,q^{D} \left| \begin{array}[t]{l} \\
L_{\mu}(\hat{\Lambda}) \end{array} \right. = \frac{q^{\frac{\langle \mu
, \mu \rangle}{2k} -\frac{\langle \Lambda
, \Lambda\rangle}{2k}}}{(q)_{\infty}^{n}} \mbox{\rm Tr}\,q^{D-
D^{\ha}} \left| \begin{array}[t]{l} \\
L_{\mu}(\hat{\Lambda})^{\ha^{+}}_{kQ} \end{array} \right.,\label{5.9}\ee
where $(q)_{\infty} := \prod_{l\geq 0}(1 - q^{l}).$ On the other hand,
one defines the string function $c^{\hat{\Lambda}}_{\mu}(q)$ as
follows (cf. [KP] or [K] Section 12.7; departing from the tradition,
we shall use subscript $\mu \in P,$ rather than $\hat{\mu} =
k\hat{\Lambda}_{0} + \mu)$:
\be c^{\hat{\Lambda}}_{\mu}(q):= q^{\frac{\langle \Lambda  +
\rho,\Lambda + \rho \rangle}{2(k+h^{\vee})} - \frac{\langle \rho ,
\rho \rangle}{2h^{\vee}} - \frac{\langle \mu , \mu \rangle}{2k}
}\;\mbox{\rm Tr}\,q^{D} \left| \begin{array}[t]{l} \\
L_{\mu}(\hat{\Lambda}) \end{array} \right.  =\label{5.10}\ee
$$= q^{\frac{\langle
\Lambda + 2 \rho,\Lambda \rangle}{2(k+h^{\vee})} -
\frac{1}{24}\frac{({\rm dim}\,\g) k}{k + h^{\vee}} - \frac{\langle \mu
, \mu \rangle}{2k} }\; \mbox{\rm Tr}\,q^{D} \left| \begin{array}[t]{l}
\\ L_{\mu}(\hat{\Lambda}) \end{array} \right.  ,$$
(the last identity follows from the strange formula of Freudenthal-de
Vries). Hence, from (\ref{5.9}),
\be c^{\hat{\Lambda}}_{\mu}(q) = \frac{q^{\frac{\langle \Lambda  +
\rho,\Lambda + \rho \rangle}{2(k+h^{\vee})} - \frac{\langle \rho ,
\rho \rangle}{2h^{\vee}} - \frac{\langle \Lambda , \Lambda \rangle}{2k}
}}{(q)_{\infty}^{n}} \;\mbox{\rm Tr}\,q^{D- D^{\ha}} \left|
\begin{array}[t]{l} \\ L_{\mu}(\hat{\Lambda})^{\ha^{+}}_{kQ}
\end{array}\right..\label{5.11}\ee\vspace{5mm}

\noindent {\it Remark 5.1} Observe that when the  level $k$ equals
$2,$ formulas
(\ref{5.11}) and (\ref{5.9}) (with the restriction (\ref{5.8}))
provide combinatorial expression for {\em every}
string function  $c^{\hat{\Lambda}}_{\mu}(q)$ corresponding to  a generic
dominant integral weight
$\hat{\Lambda} = \hat{\Lambda}_{i} +\hat{\Lambda}_{j},$ $0 \leq i,j
\leq n$: Due to the cyclic automorphism (of order $n+1)$  of the
Dynkin diagram of $\ga,$ the generic string function equals (up to a
power of $q$) a string function of the type $c^{\hat{\Lambda}_{0} +
\hat{\Lambda}_{j}}_{\mu}(q),$ $0 \leq j \leq n,$ considered
above. \vspace{5mm}

Armed with an explicit expression for the string functions, we are
only a step away from writing the character of the whole standard
module: Since $\rho(k\alpha) $ acts nontrivially only on the right
factor of the decomposition (\ref{1.4}), one concludes from
(\ref{2.201}), (\ref{2.203}) that
\be [ D , \rho(k\alpha)] \left| \begin{array}[t]{l} \\
L_{\mu}(\hat{\Lambda}) \end{array} \right. = q^{\frac{\langle
\mu + k\alpha , \mu + k\alpha\rangle}{2k} - \frac{\langle
\mu , \mu \rangle}{2k}}\rho(k\alpha) \left| \begin{array}[t]{l} \\
L_{\mu}(\hat{\Lambda}) \end{array} \right. =  \label{5.12}\ee
$$= q^{\frac{k}{2}\langle
\alpha + \frac{\mu}{k}, \alpha + \frac{\mu}{k}\rangle - \frac{\langle
\mu , \mu \rangle}{2k}}\rho(k\alpha) \left| \begin{array}[t]{l} \\
L_{\mu}(\hat{\Lambda}) \end{array} \right. = q^{\frac{k}{2}\langle
\alpha , \alpha\rangle + \langle \alpha , \mu\rangle} \rho(k\alpha)\left|
\begin{array}[t]{l} \\
L_{\mu}(\hat{\Lambda}) \end{array} \right.,  $$
where $\rho (k\alpha) :=
\underbrace{e^{\alpha}\otimes \ldots
\otimes e^{\alpha}}_{k\;factors}, \; \alpha \in Q.$  (cf.
(\ref{1.8})). Therefore (\ref{1.4}), (\ref{1.10}), (\ref{1.11}) and
(\ref{5.10})
imply that
\be \mbox{ch}\,L(\hat{\Lambda}) =  \mbox{ch}\,L(k_{0}\hat{\Lambda}_{0}
+ k_{j} \hat{\Lambda}_{j}) := \mbox{\rm
Tr}\,q^{D}\,\prod_{i=1}^{n} y_{i}^{h_{\Lambda_{i}}} \left|
\begin{array}[t]{l} \\
L(\hat{\Lambda}) \end{array} \right. = \label{5.14}\ee
$$ = \sum_{\mu  \in \Lambda + Q/kQ}  \mbox{\rm
Tr}\,q^{D} \left|
\begin{array}[t]{l} \\
L_{\mu}(\hat{\Lambda}) \end{array} \right. q^{- \frac{\langle \mu ,
\mu \rangle}{2k}}\Theta_{\mu}(q,y) = $$
$$= q^{- \frac{\langle \Lambda  +
\rho,\Lambda + \rho \rangle}{2(k+h^{\vee})} + \frac{\langle \rho ,
\rho \rangle}{2h^{\vee}}} \sum_{\mu \in \Lambda + Q/kQ}
c_{\mu}^{\hat{\Lambda}}(q)\Theta_{\mu}(q,y),$$
where $\Theta_{\mu}(q,y)$ is the classical theta function of degree $k$
(cf. e.g. [K]
Chapters 12, 13):
\be \Theta_{\mu}(q,y):=  q^{\frac{\langle \mu ,
\mu \rangle}{2k}} \sum_{\alpha \in Q} q^{\frac{k}{2} \langle \alpha ,
\alpha \rangle + \langle \alpha , \mu\rangle } \prod_{i=1}^{n}
y_{i}^{\langle \Lambda_{i}, k\alpha + \mu \rangle}  = \label{5.15}\ee
$$= \sum_{\gamma \in Q + \frac{\mu}{k}} q^{\frac{k}{2}\langle
\gamma , \gamma \rangle} \prod_{i=1}^{n} y_{i}^{k\langle \Lambda_{i},
\gamma\rangle}.$$
Formula (\ref{5.14}) is of course the familiar  expression for the
normalized character of standard module in terms of string functions
and theta functions (cf. [K] (12.7.12)). Note that the
explicit combinatorial formula for $c_{\mu}^{\hat{\Lambda}}$ is given
by (\ref{5.11}) and (\ref{5.7}) with the additional restriction
(\ref{5.8}) imposed.

   If we need only $\mbox{ch}\,L(\hat{\Lambda}),$ we can avoid any
reference to $D - D^{\ha}$-characters and use directly the
$D$-character (\ref{5.5}) of the principal subspace
$W(\hat{\dot{\Lambda}})$ (copied from [GeI] (5.27)) as well as the
above theta function which incorporates as usual the contributions of
the operators $\rho(k\alpha),\; \alpha \in Q.$ This is because by its
very definition (\ref{1.6}), (\ref{1.7}), the projection
$\pi_{L(\hat{\Lambda})^{\ha^{+}}}$ is $D$-invariant.
In other words, if we set $\dot{\mu} := \mu -
\Lambda_{j_{k}},\; \mu \in P,$ like in  (\ref{4.9})  and denote
by $Q_{(+)}$ the monoid (with 0) generated by positive roots, we have from
(\ref{1.4}), Definition 4.1, (\ref{4.10}), Theorem \ref{The
4.2}, (\ref{5.5}) and (\ref{5.12}) that
\be \mbox{ch}\,L(\hat{\Lambda})= \mbox{ch}\,L(k_{0}\hat{\Lambda}_{0}
+ k_{j} \hat{\Lambda}_{j}) := \mbox{\rm
Tr}\,q^{D}\,\prod_{i=1}^{n} y_{i}^{h_{\Lambda_{i}}} \left|
\begin{array}[t]{l} \\
L(\hat{\Lambda}) \end{array} \right. = \label{5.16}\ee
$$ = \frac{1}{(q)_{\infty}^{n}} \sum_{\mu  \in \Lambda +
Q_{(+)}}  \mbox{\rm
Tr}\,q^{D} \left|
\begin{array}[t]{l} \\
W_{\dot{\mu}}(\hat{\dot{\Lambda}}) \end{array} \right. q^{- \frac{\langle \mu ,
\mu \rangle}{2k}}\Theta_{\mu}(q,y) = $$
$$= \frac{1}{(q)_{\infty}^{n}}  \sum_{\begin{array}{c}{\scriptstyle
p_{1}^{(1)},\ldots , p_{1}^{(k-1)} \geq 0}\\{\scriptstyle \ldots \ldots
\ldots\ldots}\\{\scriptstyle p_{n}^{(1)},\ldots
, p_{n}^{(k-1)} \geq 0 }\end{array}}
\frac{q^{\frac{1}{2}\sum_{l,m =
1, \ldots, n}^{s,t = 1,\ldots , k-1}\;A_{lm} B^{st}
p_{l}^{(s)}p_{m}^{(t)}}}{\prod_{i=1}^{n} \prod_{s=1}^{k-1}
(q)_{p_{i}^{(s)}} }\;q^{\sum_{s=k_{0}+1}^{k-1} (s-k_{0})
p_{j}^{(s)}}$$ $$\sum_{\alpha \in Q} q^{\frac{k}{2}\langle \alpha ,
\alpha\rangle + \langle \alpha, \Lambda + \sum_{i=1}^{n}
r_{i}\alpha_{i}\rangle}\;\prod_{i=1}^{n} y_{i}^{\langle \Lambda_{i},
k\alpha + \Lambda + \sum_{i=1}^{n}r_{i} \alpha_{i}\rangle},$$
where as
always $r_{i} = \sum_{s=1}^{k-1}s p_{i}^{(s)}.$ Recall that
$\hat{\Lambda}$ was defined in (\ref{4.1}), $(A_{lm})_{l,m=1}^{n}$ is
the Cartan matrix of $\g$ and $B^{st} := \mbox{min}\{s,t\},$ $1\leq
s,t \leq k-1.$

   This last character formula corresponds also to a semiinfinite
monomial basis of
$L(\hat{\Lambda})$ in the spirit of Feigin and Stoyanovsky (the case
$\ga = \widehat{sl}(2,\C)$ was described by them in the announcement
[FS]). A proof of the particular case of (\ref{5.16}) for the vacuum module
($\hat{\Lambda} = k\hat{\Lambda}_{0}, $ i.e., $k_{0} = k,\; k_{j} =
0)$ was announced  in [Kir].\vspace{5mm}

\noindent {\it Example 5.1} Let $\g = sl(3,\C),\;k = 2.$ By (\ref{5.11}) and
(\ref{5.7}) and (\ref{5.8}) we have
\be c^{2\hat{\Lambda}_{0}}_{0}(q) =
\frac{q^{-\frac{2}{15}}}{(q)_{\infty}^{2}}
\sum_{\begin{array}{c}{\scriptstyle p_{1},
p_{2} \geq 0} \\ {\scriptstyle p_{1}, p_{2}\; {\rm
even}}\end{array}}\frac{q^{\frac{1}{2}(p_{1}^{2} + p_{2}^{2} -
p_{1}p_{2})}}{(q)_{p_{1}} (q)_{p_{2}}},\label{5.17}\ee
\be c^{2\hat{\Lambda}_{0}}_{\alpha_{1} + \alpha_{2}}(q) =
\frac{q^{-\frac{2}{15}}}{(q)_{\infty}^{2}}
\sum_{\begin{array}{c}{\scriptstyle p_{1},
p_{2} \geq 0} \\ {\scriptstyle p_{1}, p_{2}\; {\rm
odd}}\end{array}}\frac{q^{\frac{1}{2}(p_{1}^{2} + p_{2}^{2} -
p_{1}p_{2})}}{(q)_{p_{1}} (q)_{p_{2}}},\label{5.18}\ee
\be c^{\hat{\Lambda}_{0} + \hat{\Lambda}_{1}}_{\Lambda_{1}}(q) =
\frac{q^{-\frac{1}{30}}}{(q)_{\infty}^{2}}
\sum_{\begin{array}{c}{\scriptstyle p_{1},
p_{2} \geq 0} \\ {\scriptstyle p_{1},p_{2}\; {\rm
even}} \end{array}}\frac{q^{\frac{1}{2}(p_{1}^{2} + p_{2}^{2} -
p_{1}p_{2}- p_{1})}}{(q)_{p_{1}} (q)_{p_{2}}},\label{5.19}\ee
\be c^{\hat{\Lambda}_{0} + \hat{\Lambda}_{1}}_{\Lambda_{1} +\alpha_{2}}(q) =
\frac{q^{-\frac{1}{30}}}{(q)_{\infty}^{2}}
\sum_{\begin{array}{c}{\scriptstyle p_{1},
p_{2} \geq 0} \\ {\scriptstyle p_{1}\;(p_{2})\; {\rm even}\;({\rm
odd})}\end{array}}\frac{q^{\frac{1}{2}(p_{1}^{2} + p_{2}^{2} -
p_{1}p_{2}- p_{1})}}{(q)_{p_{1}} (q)_{p_{2}}},\label{5.20}\ee Due to
symmetries, the above string functions determine all the other string
functions at level $2$ (cf.  for example [KP] Section 4.6; using the
expressions of Kac and Peterson, the above formulas were verified on
Maple up to ${\cal O}(50)).$ In particular, due to the cyclic automorphism
of the affine Dynkin diagram, one has $c^{\hat{\Lambda}_{0} +
\hat{\Lambda}_{1}}_{\Lambda_{1}}(q) = c^{\hat{\Lambda}_{0} +
\hat{\Lambda}_{2}}_{\Lambda_{2}}(q) = c^{\hat{\Lambda}_{1} +
\hat{\Lambda}_{2}}_{\Lambda_{1} + \Lambda_{2}}(q).$ But
notice that $\hat{\Lambda}_{1} + \hat{\Lambda}_{2}$ is not among the
highest weights of type (\ref{4.1}) for which our basis theorems and
character formulas hold (cf. Remark 5.1). Thanks to the simplicity of
the case, one can nevertheless  find easily a combinatorial basis and
write down the corresponding character formula:
\be  \mbox{Tr}\,q^{D- D^{\ha}}\left | \begin{array}[t]{l} \\
L(\hat{\Lambda}_{1} + \hat{\Lambda}_{2})^{\ha^{+}}_{2Q}
\end{array} \right. = \label{5.21}\ee
$$= \sum_{p_{1},
p_{2} \geq 0}\frac{q^{\frac{1}{2}(p_{1}^{2} + p_{2}^{2} -
p_{1}p_{2})}}{(q)_{p_{1}} (q)_{p_{2}}} q^{p_{1}- \frac{1}{2}(p_{1} + p_{2})}
 + \sum_{p_{1},
p_{2} \geq 0}\frac{q^{\frac{1}{2}[(p_{1} + 1)^{2} + p_{2}^{2} -
(p_{1} + 1)p_{2}]}}{(q)_{p_{1}} (q)_{p_{2}}} q^{p_{2} -
\frac{1}{2}[(p_{1} + 1) + p_{2}]}.$$
The first term counts only monomials which do not include a
quasi-particle  $\pi_{L(\hat{\Lambda}_{1} +
\hat{\Lambda}_{2})^{\ha^{+}}_{2Q}} \cdot x_{\alpha_{1}}(-1),$  while
the second term counts only monomials which contain such a
quasi-particle.

\section {Appendix}

$\begin{array}{|c|c|c|} \hline \mbox{ \em color-} & \mbox{ \em energy}
&
\mbox{ \em basis}
\\ \mbox{\em -type} & & \\ \hline\hline   (1;2) & 3/2 & (1_{\alpha_{2}}
-3_{\alpha_{1}}-1_{\alpha_{1}}) \\ \cline{2-3} & 5/2 &
(1_{\alpha_{2}}-4_{\alpha_{1}}-1_{\alpha_{1}}), (0_{\alpha_{2}}-
3_{\alpha_{1}}-1_{\alpha_{1}})\\ \cline{2-3} & 7/2 &
(1_{\alpha_{2}}-5_{\alpha_{1}}-1_{\alpha_{1}}),
(1_{\alpha_{2}}-4_{\alpha_{1}}-2_{\alpha_{1}}),
(0_{\alpha_{2}}-4_{\alpha_{1}}-1_{\alpha_{1}}),\\ & &
(1_{\alpha_{2}}-3_{\alpha_{1}}-1_{\alpha_{1}}) \\ \hline (2;2) & 2 &
(-1_{\alpha_{2}}1_{\alpha_{2}}-3_{\alpha_{1}}-1_{\alpha_{1}}) \\
\cline{2-3} & 3 &
(-2_{\alpha_{2}}1_{\alpha_{2}}-3_{\alpha_{1}}-1_{\alpha_{1}}),
(-1_{\alpha_{2}} 1_{\alpha_{2}} -4_{\alpha_{1}} -1_{\alpha_{1}})
\\ \cline{2-3} & 4 &  (-3_{\alpha_{2}}1_{\alpha_{2}}-3_{\alpha_{1}}
-1_{\alpha_{1}}),
(-2_{\alpha_{2}}1_{\alpha_{2}}-4_{\alpha_{1}}-1_{\alpha_{1}}) ,
(-2_{\alpha_{2}}0_{\alpha_{2}}-3_{\alpha_{1}} -1_{\alpha_{1}}), \\ & &
(-1_{\alpha_{2}}1_{\alpha_{2}} -5_{\alpha_{1}} -1_{\alpha_{1}}),
(-1_{\alpha_{2}}1_{\alpha_{2}}-4_{\alpha_{1}}-2_{\alpha_{1}}) \\
\hline \end{array}$ \begin {center} Table 1 \end{center}\vspace{5mm}

\noindent$\begin{array}{|c|c|c|c|} \hline  \mbox{ \em color-} &  \mbox{ \em
energy} & \mbox{\em color-} & \mbox{\em basis}\\ \mbox{\em -type} & &
\mbox{\em -charge-} & \\ & & \mbox{\em -type} &\\
\hline \hline (1;2) & 1 &(1;2)& (0_{\alpha_{2}}-2_{2\alpha_{1}})\\
\cline{2-4} & 2
& (1;1,1)& (1_{\alpha_{2}}-3_{\alpha_{1}}-1_{\alpha_{1}})\\
\cline{3-4}& & (1;2)&
(0_{\alpha_{2}}-3_{2\alpha_{1}}),(-1_{\alpha_{2}}-2_{2\alpha_{1}})\\
\cline{2-4} & 3 &(1;1,1) &  (1_{\alpha_{2}}-4_{\alpha_{1}}-1_{\alpha_{1}}),
(0_{\alpha_{2}}-3_{\alpha_{1}}-1_{\alpha_{1}})\\\cline{3-4} & & (1;2)
& (0_{\alpha_{2}}-4_{2\alpha_{1}}), (-1_{\alpha_{2}}
-3_{2\alpha_{1}}), (-2_{\alpha_{2}}-2_{2\alpha_{1}}) \\ \cline{2-4} &
4 & (1;1,1) & (1_{\alpha_{2}} -4_{\alpha_{1}} -2_{\alpha_{1}}),
(1_{\alpha_{2}} -5_{\alpha_{1}} -1_{\alpha_{1}}),\\ & &
&(0_{\alpha_{2}} -4_{\alpha_{1}}-1_{\alpha_{1}}),
(-1_{\alpha_{2}}-3_{\alpha_{1}}-1_{\alpha_{1}})\\\cline{3-4}& & (1;2)&
(0_{\alpha_{2}} -5_{2\alpha_{1}}), (-1_{\alpha_{2}} -4_{2\alpha_{1}}),
(-2_{\alpha_{2}} -3_{2\alpha_{1}}),\\ & & & (-3_{\alpha_{2}}
-2_{2\alpha_{1}})\\
\hline (2;2) & 2/3 & (2;2) &
(0_{2\alpha_{2}}-2_{2\alpha_{1}})\\\cline{2-4}& 5/3
& (2;2) & (0_{2\alpha_{2}} -3_{2\alpha_{1}}), (-1_{2\alpha_{2}}
-2_{2\alpha_{1}})\\ \cline{2-4} & 8/3 & (1,1;1,1) & (-1_{\alpha_{2}}
1_{\alpha_{2}}-3_{\alpha_{1}}-1_{\alpha_{1}})\\ \cline{3-4} & &
(2;1,1) & (0_{2\alpha_{2}}-3_{\alpha_{1}}-1_{\alpha_{1}})\\
\cline{3-4} & & (1,1;2) & (-2_{\alpha_{2}} 0_{\alpha_{2}}
-2_{2\alpha_{1}})\\ \cline{3-4} & & (2;2) & (0_{2\alpha_{2}}
-4_{2\alpha_{1}}), (-1_{2\alpha_{2}} -3_{2\alpha_{1}}),
(-2_{2\alpha_{2}} -2_{2\alpha_{1}})\\ \cline{2-4} & 11/3 & (1,1;1,1) &
(-1_{\alpha_{2}} 1_{\alpha_{2}} -4_{\alpha_{1}}
-1_{\alpha_{1}}),(-2_{\alpha_{2}} 1_{\alpha_{2}}
-3_{\alpha_{1}}-1_{\alpha_{1}})\\ \cline{3-4} & & (2;1,1) &
(0_{2\alpha_{2}} -4_{\alpha_{1}}-1_{\alpha_{1}}), (-1_{2\alpha_{2}}
-3_{\alpha_{1}} -1_{\alpha_{1}}) \\ \cline{3-4} & & (1,1;2) &
(-2_{\alpha_{2}} 0_{\alpha_{2}} -3_{2\alpha_{1}}),(-3_{\alpha_{2}}
0_{\alpha_{2}} -2_{2\alpha_{1}}) \\ \cline{3-4} & & (2;2) &
(0_{2\alpha_{2}} -5_{2\alpha_{1}}), (-1_{2\alpha_{2}}
-4_{2\alpha_{1}}), (-2_{2\alpha_{2}} -3_{2\alpha_{1}}),\\ & & &
(-3_{2\alpha_{2}} -2_{2\alpha_{1}})\\ \hline \end{array}$
\begin {center} Table 2 \end{center}\vspace{10mm}
{\sc Dept. of Mathematics, Rutgers University, New Brunswick, NJ
08903}\\georgiev@math.rutgers.edu


\vspace{10mm}




\begin{thebibliography}{AAAAAA}
\bibitem[A]{1} {\sc G.E. Andrews,} {\it The theory of partitions},
Addison-Wesley, 1976.
\bibitem[ANOT]{} {\sc T. Arakawa, T. Nakanishi, K. Ooshima and A.
Tsuchiya,} talk in the meeting of JMS at Ritsumeikan University (1995).
\bibitem[BG]{} {\sc E. Baver and D. Gepner,} {\it Fermionic sum
representations for the Virasoro characters of the unitary
superconformal unitary models,} (hep-th/9502118).
\bibitem[BLS1]{} {\sc P. Bouwknegt, A. Ludwig and K. Schoutens,} {\it
Spinon basis
for higher level $SU(2)$ WZW models,} (hep-th/9412108); {\it Spinon
bases, Yangian symmetry and fermionic representations of Virasoro
characters in conformal field theory,} Phys. Lett.  {\bf 338B} (1994),
448; (hep-th/9406020); {\it Affine and Yangian symmetries in
$SU(2)_{1}$ conformal field theory,} (hep-th/9412199);
\bibitem[BLS2]{} {\sc P. Bouwknegt, A. Ludwig and K. Schoutens,} {\it
Spinon basis for $(\widehat{sl_{2}})_{k}$ integrable highest weight
modules and new character formulas,} (hep-th/9504074).
\bibitem[BM]{} {\sc A. Berkovich and B. McCoy,} {\it  Continued fractions and
fermionic representations for characters of $M(p,p')$ minimal models,}
(hep-th/9412030).
\bibitem[BPS]{} {\sc D. Bernard, V.   Pasquier   and D. Serban,} {\it
Spinons in
Conformal Field Theory,} Nucl. Phys.  {\bf B428} (1994), 612,
(hep-th/9404050).
\bibitem[DKKMM]{} {\sc S. Dasmahapatra, R.  Kedem,  T.R.   Klassen,
B.  McCoy and E. Melzer,} {\it Quasi-particles, conformal field theory
and q series,} Int. J. Mod. Phys. {\bf B7} (1993), 3617,
(hep-th/9303013); {\it Virasoro characters from Bethe equations for
the critical ferromagnetic three-state Potts model,} J. Stat. Phys.
{\bf 74} (1994), 239, (hep-th/9304150).
\bibitem[DL]{} {\sc C. Dong and J. Lepowsky,} {\it Generalized Vertex
Algebras and Relative Vertex Operators}, Progress in Math., Vol.112,
Birkhauser, Boston, 1993.
\bibitem[FIJKMY]{} {\sc O. Foda, K. Iohara, M. Jimbo, R. Kedem, T.
Miwa and H. Yan,} {\it Notes on highest weight modules of the elliptic
algebra $A_{q,p}(\widehat{sl}(2)),$} (hep-th/9405058).
\bibitem[FK]{18} {\sc I. Frenkel and V. Kac,} {\it  Basic representations of
affine Lie algebras and dual resonance models,} Invent. Math.
{\bf 62} (1880), 23-66.
\bibitem[FLM]{6} {\sc I. Frenkel, J.  Lepowsky  and A. Meurman,} {\it  Vertex
Operator Algebras and the Monster,}  Pure and Appl. Math. {\bf 134},
Academic Press 1988.
\bibitem[FS]{} {\sc B.L. Feigin and A.V. Stoyanovsky,} {\it
Quasi-particles models for the representations of Lie algebras and
geometry of flag manifold,} (hep-th/9308079); cf.
also the short version: {\it Functional models for representations of
current algebras and semi-infinite Schubert cells,}  Funct. Anal.
Appl. {\bf 28}  No.1 (1994), 55.
\bibitem[FW]{} {\sc O. Foda   and S. O. Warnaar,} {\it  A bijection
which implies
Melzer's polynomial identities: the $\chi_{1,1}^{(p,p+1)}$ case,}
(hep-th/9501088).
\bibitem[G]{9} {\sc D. Gepner, } {\it New conformal field theories associated
with Lie algebras and their partition functions,} Nucl. Phys. B
{\bf 290} (1987), 10-24.
\bibitem[GeI]{}  {\sc G. Georgiev,} {\it  Parafermionic constructions of
modules for  infinite - dimensional  Lie algebras, I.  Principal
subspace,} (hep-th/9412054).
\bibitem[H]{} {\sc F.D.M. Haldane,} {\it  ''Fractional statistics'' in
arbitrary
dimensions: a generalization of the Pauli principle,} Phys. Rev.
Lett.  {\bf 67} (1991), 937-940.
\bibitem[I]{} {\sc S. Iso,} {\it Anyon basis of $c = 1$ conformal field
theory,} (hep-th/9411051).
\bibitem[K]{10} {\sc V.G. Kac,}  {\it Infinite-dimensional Lie
algebras,} Cambridge University Press, Cambridge, 1990.
\bibitem[Kir]{} {\sc A.N. Kirillov,} {\it Dilogarithm identities,} (
hep-th/9408113).
\bibitem[KKMM]{3} {\sc R. Kedem, T.    Klassen,
 B. McCoy and E. Melzer,} {\it Fermionic quasi-particle
representations for characters of $\frac{(G^{(1)})_{1} \times
(G^{(1)})_{1}}{(G^{(1)})_{2}},$} Phys. Lett.  {\bf B304} (1993),
263-270, (hep-th/9211102); {\it Fermionic sum representations for
conformal field theory characters,} Phys. Lett. {\bf B307} (1993),
68-76, (hep-th/9301046).
\bibitem[KMM]{} {\sc R. Kedem, B.   McCoy  and E.  Melzer,} {\it   The Sums of
Rogers, Schur
and Ramanujan and the Bose-Fermi correspondence in 1+1-dimensional
quantum field theory,} (hep-th/9304056).
\bibitem[KNS]{11} {\sc A. Kuniba, T. Nakanishi and J. Suzuki,} {\it
Characters of conformal field theories from thermodynamic Bethe
Ansatz,}  Mod. Phys. Lett.  {\bf A8} (1993), 1649-1660, (hep-th/9301018).
\bibitem[KP]{12} {\sc V. Kac and D. Peterson,} {\it Infinite-dimensional
Lie algebras, theta functions and modular forms,}  Adv. Math.  {\bf
53} (1984), 125-264.
\bibitem[LP]{13} {\sc J. Lepowsky and M. Primc,} {\it  Structure of the
standard modules for the affine algebra $A_{1}^{(1)},$} Contemp.
Math. {\bf 46}, 1985.
\bibitem[LW]{56} {\sc J. Lepowsky  and  R. Wilson,} {\it   A new family
of algebras underlying the Rogers-Ramanujan identities and
generalizations,} Proc.  Natl. Acad. Sci. USA {\bf 78} (1981),
7245-7248; {\it   The
structure of standard
modules, I: Universal algebras and Rogers-Ramanujan identities,}
Invent. Math.  {\bf 77} (1984), 199-290; {\it   The
structure of standard
modules, II: The case $A_{1}^{(1)},$ principal gradation,} Invent.  Math.
{\bf 79} (1985), 417-442.
\bibitem[M]{} {\sc E. Melzer,} {\it   Fermionic character sums and
the corner
transfer matrix,} Int. J. Mod. Phys.  {\bf A9} (1994), 1115-1136,
(hep-th/9305114).
\bibitem[NY]{} {\sc A. Nakayashiki and Y. Yamada,} {\it Crystalizing
the spinon basis,} (hep-th/9504052).
\bibitem[P]{55} {\sc M. Primc,} {\it  Vertex operator construction of standard
modules for $A_{n}^{(1)},$} Pacific J. Math. {\bf 162}  No.1
(1994), 143-187.
\bibitem[S]{19} {\sc G. Segal,} {\it  Unitary representations of some
infinite-dimensional groups,} Commun. Math. Phys.  {\bf 80} (1981), 301.
\bibitem[W]{} {\sc S. O. Warnaar,} {\it Fermionic solutions of the
Andrews-Baxter-Forrester model I: unification of TBA and CTM methods,}
(hep-th/9501134).
\bibitem[WP]{} {\sc S. O. Warnaar and P. Pearce,} {\it  A-D-E polynomial and
Rogers - Ramanujan identities,} (hep-th/9411009).
\bibitem[ZF]{16} {\sc  A.B. Zamolodchikov and V.A. Fateev,} {\it
Non-local (parafermion) currents in two-dimensional quantum field
theory and self-dual critical points in $\Z_{n}$- symmetric classical
systems,}   Sov. Phys. JETP  {\bf 62} (1985), 215.


\end{thebibliography}
\end{document}